# ON THE DISCONNECTION OF A DISCRETE CYLINDER BY A BIASED RANDOM WALK

By David Windisch

*ETH Zürich*

We consider a random walk on the discrete cylinder $(\mathbb{Z}/N\mathbb{Z})^d \times \mathbb{Z}$, $d \geq 3$ with drift $N^{-d\alpha}$ in the $\mathbb{Z}$-direction and investigate the large $N$-behavior of the disconnection time $T_N^{\mathrm{disc}}$, defined as the first time when the trajectory of the random walk disconnects the cylinder into two infinite components. We prove that, as long as the drift exponent $\alpha$ is strictly greater than 1, the asymptotic behavior of $T_N^{\mathrm{disc}}$ remains $N^{2d+o(1)}$, as in the unbiased case considered by Dembo and Sznitman, whereas for $\alpha < 1$, the asymptotic behavior of $T_N^{\mathrm{disc}}$ becomes exponential in $N$.

**1. Introduction.** Informally, the object of our study can be described as follows: a particle feeling a drift moves randomly through a cylindrical object, and damages every visited point. How long does it take until the cylinder breaks apart, and how does the answer to this question depend on the drift felt by the particle? This is a variation on the problem of "the termite in a wooden beam" considered by Dembo and Sznitman [4].

We henceforth consider the *discrete cylinder*

$$(1.1) \qquad E = \mathbb{T}_N^d \times \mathbb{Z}, \qquad d \geq 1,$$

where $\mathbb{T}_N^d$ denotes the $d$-dimensional integer torus $\mathbb{T}_N^d = (\mathbb{Z}/N\mathbb{Z})^d$. The disconnection time of the cylinder $E$ by a simple (unbiased) random walk was introduced by Dembo and Sznitman in [4], where it was shown that its asymptotic behavior is approximately $N^{2d} = |\mathbb{T}_N^d|^2$ as $N \to \infty$ when $d \geq 1$. This result was extended by Sznitman in [12] to a wide class of bases of $E$ with uniformly bounded degree as $N \to \infty$. Similar models related to interfaces created by simple random walk trajectories have been studied by Benjamini and Sznitman [3] and Sznitman [13]. The former of these two works has led Dembo and Sznitman [5] to sharpen their lower bound on the

---

disconnection time of $E$ for large $d$. Here we investigate the disconnection time for a random walk with bias into the $\mathbb{Z}$-direction.

We now proceed to the precise description of the problem studied in the present work. The cylinder $E$ is equipped with the Euclidean distance $|\cdot|$ and the natural product graph structure, for which all vertices $x_1$, $x_2 \in E$ with $|x_1 - x_2| = 1$ are connected by an edge. The (discrete-time) random walk with drift $\Delta \in [0, 1)$ is the Markov chain $(X_n)_{n \geq 0}$ on $E$ with starting point $x \in E$ and transition probability

$$(1.2) \quad p_X(x_1, x_2) = \frac{1 + \Delta(\pi_{\mathbb{Z}}(x_2 - x_1))}{2d + 2} \mathbf{1}_{\{|x_1 - x_2| = 1\}}, \qquad x_1, x_2 \in E,$$

where $\pi_{\mathbb{Z}}$ denotes the projection from $E$ onto $\mathbb{Z}$. The process is defined on a suitable filtered probability space $(\Omega_N, (\mathcal{F}_n)_{n \geq 0}, P_x^\Delta)$ (see Section 2 for details). In particular, under $P_0^0$, $X$ is the ordinary simple random walk on $E$. We say that a set $K \subseteq E$ *disconnects* $E$ if $\mathbb{T}_N^d \times (-\infty, -M]$ and $\mathbb{T}_N^d \times [M, \infty)$ are contained in two distinct components of $E \setminus K$ for large $M \geq 1$. The central object of interest is the *disconnection time*

$$(1.3) \qquad\qquad T_N^{\mathrm{disc}} = \inf\{n \geq 0 : X([0, n]) \text{ disconnects } E\}.$$

We consider drifts of the form $N^{-d\alpha} = |\mathbb{T}_N^d|^{-\alpha}$, $\alpha > 0$. Our main result shows that the asymptotic behavior of $T_N^{\mathrm{disc}}$ as $N \to \infty$ is the same as in the case without drift considered in [4] as long as $\alpha > 1$, and becomes exponential in $N$ when $\alpha < 1$:

THEOREM 1.1 ($d \geq 3$, $\alpha > 0$, $\varepsilon > 0$).

$$(1.4) \quad \begin{aligned} &\text{For } \alpha > 1, &&N^{2d-\varepsilon} \leq T_N^{\mathrm{disc}} \leq N^{2d+\varepsilon}, \\ &\text{for } \alpha < 1, &&\exp\{N^{d(1-\alpha-\varphi(\alpha))-\varepsilon}\} \leq T_N^{\mathrm{disc}} \leq \exp\{N^{d(1-\alpha)+\varepsilon}\}, \end{aligned}$$

*with probability tending to 1 as $N \to \infty$, where the continuous function $\varphi : (0, 1) \to (0, \frac{1}{d-1})$ is defined by*

$$(1.5) \quad \begin{aligned} \varphi(\alpha) = {}&\alpha \mathbf{1}_{\{0 < \alpha < \alpha_*\}} + \left(\frac{1}{d} + \frac{\alpha}{d-1} - \alpha\right) \mathbf{1}_{\{\alpha_* < \alpha < 1/d\}} \\ &+ \frac{1-\alpha}{(d-1)^2} \mathbf{1}_{\{1/d \leq \alpha < 1\}}, \end{aligned}$$

*for $\alpha_* = \frac{1}{d(2-1/(d-1))}$. In particular, $\varphi$ satisfies $\lim_{\alpha \to 0} \varphi(\alpha) = \lim_{\alpha \to 1} \varphi(\alpha) = 0$ [see Figure 1 for an illustration of the region between $1 - \alpha - \varphi(\alpha)$ and $1 - \alpha$].*

We now outline the ideas entering the proof of this result. The upper bounds on $T_N^{\mathrm{disc}}$ are derived in Theorem 3.1. The proof of this theorem is based on the simple observation that the cylinder $E$ is disconnected as soon as a slice of the form $\mathbb{T}_N^d \times \{z\} \subseteq E$ is completely covered by the walk. We thus show that the trajectory of the random walk $X$ up to time $N^{2d+\varepsilon}$ (for $\alpha > 1$), respectively $\exp\{N^{d(1-\alpha)+\varepsilon}\}$ (for $\alpha < 1$), does cover such a slice with probability tending to 1 as $N \to \infty$. To this end, we fix the slice $\mathbb{T}_N^d \times \{0\}$ and record visits made to it by $X$, where we only count visits with a sufficient time of "relaxation" in between. The process recording these visits is defined as $(V, 0)$ [cf. (3.8)]. Once we have checked that $V$ forms a Markov chain on $\mathbb{T}_N^d$ in Lemma 3.5, we can infer from the coupon-collector-type estimate (3.9) on the cover time that after a certain "critical" number of visits, the slice $\mathbb{T}_N^d \times \{0\}$ is covered with overwhelming probability by $(V, 0)$, hence by $X$. Since the same estimates apply to any slice $\mathbb{T}_N^d \times \{z\}$, $z \in \mathbb{Z}$, we are left with the one-dimensional problem of finding an upper bound on the time until sufficiently many such visits occur for some slice $\mathbb{T}_N^d \times \{z\}$.

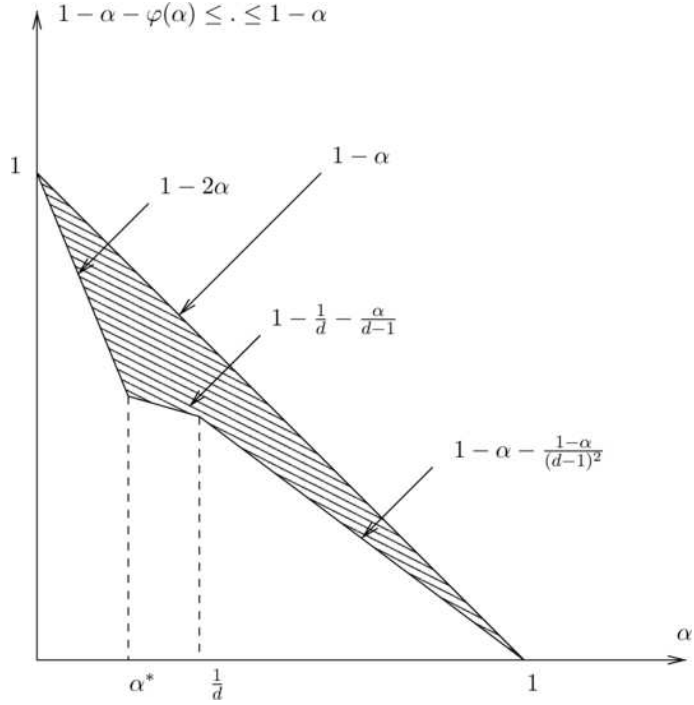

Fig. 1. *The shaded region lies between the exponents of the upper and lower bounds in Theorem 1.1 for $\alpha \in (0,1)$.*

Let us now describe the ideas involved in the more delicate derivation of the lower bounds. In this work, we reduce the problem of finding a lower bound on $T_N^{\mathrm{disc}}$ to a large deviations problem concerning the disconnection of a certain finite subset of $E$ by excursions of an unbiased simple random walk, and then derive estimates on this large deviations problem. Let us describe this last problem and the reduction step in more detail. For any subsets $K$, $B \subseteq E$, $B$ finite, and $\kappa \in (0, \frac{1}{2})$, we say that $K$ *$\kappa$-disconnects* $B$ if $K$ contains the relative boundary in $B$ of a subset of $B$ with relative volume between $\kappa$ and $1 - \kappa$, that is, if there is a subset $I$ of $B$ (generally not unique) such that

$$(1.6) \qquad \kappa|B| \leq |I| \leq (1 - \kappa)|B| \quad \text{and} \quad \partial_B(I) \subseteq K,$$

where, for sets $A$, $B \subseteq E$, $|A|$ denotes the number of points in $A$ and $\partial_B(A)$ the $B$-relative boundary of $A$, that is, the set of points in $B \setminus A$ with neighbors in $A$. The set whose disconnection concerns us is

$$(1.7) \qquad B(\alpha) = \left[-\left[\frac{N}{4}\right], \left[\frac{N}{4}\right]\right]^d \times \left[-\left[\frac{N^{d\alpha \wedge 1}}{4}\right], \left[\frac{N^{d\alpha \wedge 1}}{4}\right]\right].$$

Note that in the case $\alpha \geq \frac{1}{d}$, $B(\alpha)$ becomes $B_\infty(0, [N/4])$, the closed ball of radius $[N/4]$ with respect to the $l_\infty$-distance, centered at 0. We define $U_{B(\alpha)}$ as the first time when the trajectory of the random walk $\frac{1}{3}$-disconnects $B(\alpha)$, that is,

$$(1.8) \qquad U_{B(\alpha)} = \inf\{n \geq 0 : X([0, n]) \ \frac{1}{3}\text{-disconnects } B(\alpha)\}.$$

The random walk excursions featuring in the large deviations problem are excursions in and out of slices of the form

$$(1.9) \qquad S_r = \mathbb{T}_N^d \times [-[r], [r]] \subseteq E \qquad (r > 0).$$

Finally, the crucial reduction step comes in the following theorem, proved in Section 4:

THEOREM 1.2 $(d \geq 2, \ \alpha > 0, \ \beta > 0)$. *Suppose that $f$ is a nonnegative function on $(0, \infty)^2$ such that, for $(\mathcal{R}_n)_{n \geq 1}$, $(\mathcal{D}_n)_{n \geq 1}$, the successive returns to $S_{2[N^{d\alpha \wedge 1}]}$ and departures from $S_{4[N^{d\alpha \wedge 1}]}$ [cf. (2.24)] and the stopping time defined in (1.8), one has*

$$(1.10) \qquad \varlimsup_{N \to \infty} \frac{1}{N^\xi} \log \sup_{x \in S_{2[N^{d\alpha \wedge 1}]}} P_x^0[U_{B(\alpha)} \leq \mathcal{D}_{[N^\beta]}] < 0$$

$$\text{for any } 0 < \xi < f(\alpha, \beta).$$

*If $f(\alpha, \beta) > 0$ for all $\alpha > 1$, $\beta \in (0, d - 1)$, then it follows that*

$$(1.11) \qquad P_0^{N^{-d\alpha}}[N^{2d-\varepsilon} \leq T_N^{\mathrm{disc}}] \overset{N \to \infty}{\longrightarrow} 1 \qquad \text{for any } \alpha > 1, \varepsilon > 0,$$

*while for any $f \geq 0$,*

$$(1.12) \qquad P_0^{N^{-d\alpha}}[\exp\{N^{\zeta - \varepsilon}\} \leq T_N^{\mathrm{disc}}] \overset{N \to \infty}{\longrightarrow} 1 \qquad \text{for any } \alpha > 0, \varepsilon > 0,$$

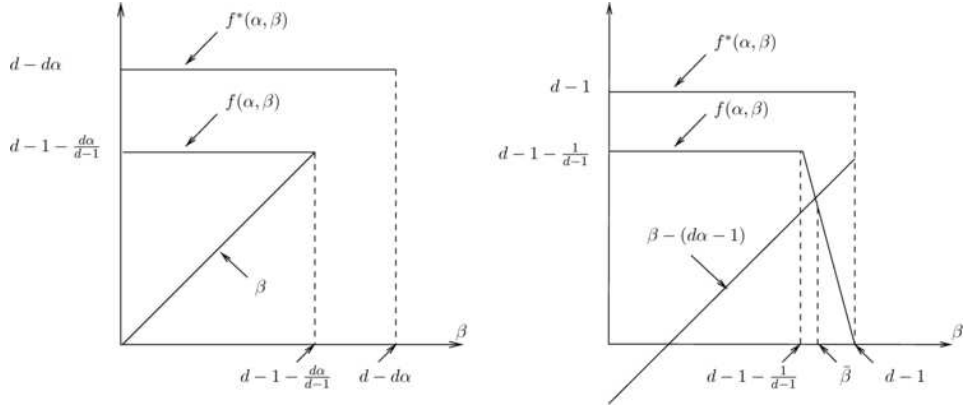

FIG. 2. *The function* $f$ *provided by Theorem* 6.1, *case* $\alpha \in (0, \frac{1}{d})$ *on the left, case* $\alpha \in [\frac{1}{d}, \infty)$ *on the right.*

*where*

$$(1.13) \qquad \zeta = \sup_{\beta > 0} g_\alpha(\beta) \quad and \quad g_\alpha(\beta) = (\beta - (d\alpha - 1)_+) \wedge f(\alpha, \beta).$$

In order to apply Theorem 1.2, one has to find a suitable nonnegative function $f$ satisfying the fundamental large deviations estimate (1.10). We show in Theorem 6.1 that (1.10) holds for the function $f$ illustrated in Figure 2. With this function $f$, the lower bound exponents $\zeta$ [in (1.12)] and $d(1 - \alpha - \varphi(\alpha))$ [in (1.4)] are related via

$$(1.14) \qquad d(1 - \alpha - \varphi(\alpha)) = \zeta \vee (d(1 - 2\alpha)\mathbf{1}_{\{\alpha < 1/d\}}),$$

as is shown in Corollary 6.3. The fact that the lower bound on $T_N^{\mathrm{disc}}$ holds with the expression $d(1 - 2\alpha)\mathbf{1}_{\{\alpha < 1/d\}}$ in (1.14) follows from the rather straightforward lower bound derived in Proposition 6.2.

We now sketch some of the techniques involved in the proof of Theorem 1.2 and the subsequent derivation of the large deviations estimate (1.10). The first step in the proof of Theorem 1.2 is a purely geometric argument in the spirit of Dembo and Sznitman [4] showing that any trajectory disconnecting $E$ must $\frac{1}{3}$-disconnect a set of the form $x_* + B(\alpha)$ (see Lemma 4.1). On the event that the walk performs no more than $[N^\beta]$ excursions between $x_* + S_{2[N^{d\alpha-1}]}$ and $x_* + S_{4[N^{d\alpha-1}]}^c$ for any $x_* \in E$ until some time $t_N$, disconnection before time $t_N$ can only occur if these at most $[N^\beta]$ excursions $\frac{1}{3}$-disconnect $x_* + B(\alpha)$ for some $x_* \in E$. One can thus apply the assumed large deviations estimate (1.10) after getting rid of the drift with the help of a Girsanov-type control (see Lemma 2.1) and applying translation invariance. It then remains to bound the probability that more than $[N^\beta]$ of the above-mentioned excursions occur for some $x_* \in E$. This can be achieved with standard estimates on one-dimensional random walk.

In order to derive the fundamental large deviations estimate (1.10), we begin with some more geometric lemmas. We show in Lemmas 5.1–5.3 that when $0 < \gamma < \gamma' < 1$, for large $N$ and any set $K$ $\frac{1}{3}$-disconnecting $B_\infty(0, [N/4])$ [cf. (1.7) and thereafter], one can find a subcube of $B_\infty(0, [N/4])$ with size $L = [N^{\gamma'}]$, so that $K$ contains a "well-spread" set of points in each of a "well-spread" collection of sub-subcubes with size $l = [N^\gamma]$ (we refer to Lemma 5.3 for the precise statement). A key ingredient for the proof of this geometric result, similar to Lemma 2.5 of [4], is an isoperimetric inequality from [6] (see Lemma 5.2). A small modification of the argument shows a similar result for $B(\alpha)$, $\alpha < \frac{1}{d}$ (see Lemma 5.4). As a consequence of these geometric results, one finds that the event under consideration in the large deviations estimate (1.10) is included in the event that the trajectory left by the $[N^\beta]$ excursions has substantial presence in many small subcubes of $B(\alpha)$. The key control on an event of this form is provided by Lemma 6.5. The main part of the argument there is to obtain a tail estimate on the number of points contained in the projection on one of the $d$-dimensional hyperplanes of the small subcubes intersected with the trajectory of the random walk stopped when exiting a large set. It follows from Khaśminskii's lemma that this number of points, divided by its expectation, is a random variable whose exponential moment is uniformly bounded with $N$. In order to bound the expected number of visited points, we use standard estimates on the Green function of the simple random walk.

An obvious question arising from Theorem 1.1 is whether one can prove the same result with $\varphi \equiv 0$ in (1.4). With Theorem 1.2, it is readily seen that this would follow if one could show that the large deviations estimate (1.10) holds with

$$(1.15) \qquad f^*(\alpha, \beta) = \begin{cases} d - (d\alpha \wedge 1), & \beta < d - (d\alpha \wedge 1), \\ 0, & \beta \geq d - (d\alpha \wedge 1); \end{cases}$$

see Figure 2. In fact, the above function $f^*$ can be shown to be the correct exponent associated to a large deviations problem similar to (1.10), where one replaces the time $U_{B(\alpha)}$ by $\mathcal{U}$, defined as the first time when the trajectory of $X$ covers $\mathbb{T}_N^d \times \{0\}$. Plainly one has $U_{B(\alpha)} \leq \mathcal{U}$, and it follows that any function $f$ in (1.10) satisfies $f(\alpha, \beta) \leq f^*(\alpha, \beta)$ for all points $(\alpha, \beta)$ of continuity of $f$; we refer to Remark 6.7 for more details. The crucial open question is therefore: are these two problems sufficiently similar for (1.10) to hold with $f^*$?

*Organization of the article.* In Section 2, we provide the definitions and the notation to be used throughout this article and prove a Girsanov-type estimate to be frequently used later on.

In Section 3, we derive the upper bounds on $T_N^{\mathrm{disc}}$ of Theorem 1.1.

In Section 4, we prove Theorem 1.2, thus reducing the derivation of a lower bound on $T_N^{\mathrm{disc}}$ to a large deviations estimate.

In Section 5, we prove several geometric lemmas in preparation of our derivation of the latter estimate.

In Section 6, we supply the key large deviations estimate in Theorem 6.1 and derive a simple lower bound on $T_N^{\mathrm{disc}}$ for large drifts. As we show, this yields the lower bounds on $T_N^{\mathrm{disc}}$ in Theorem 1.1.

*Constants.* Finally, we use the following convention concerning constants: Throughout the text, $c$ or $c'$ denote positive constants which only depend on the base-dimension $d$, with values changing from place to place. The numbered constants $c_0, c_1, \ldots$ are fixed and refer to their first place of appearance in the text. Dependence of constants on parameters other than $d$ appears in the notation. For example, $c(\gamma, \gamma')$ denotes a positive constant depending on $d$, $\gamma$ and $\gamma'$.

## 2. Definitions, notation and a useful estimate.
The purpose of this section is to set up the notation and the definitions to be used in this article and to provide a Girsanov-type estimate comparing the random walks with drift and without drift, to be frequently applied later on.

Throughout this article, we denote, for $s, t \in \mathbb{R}$, by $s \wedge t$ the minimum of $s$ and $t$, by $s \vee t$ the maximum of $s$ and $t$, by $[s]$ the largest integer satisfying $[s] \leq s$ and we set $t_+ = t \vee 0$ and $t_- = -(t \wedge 0)$. Recall that we introduced the cylinder $E$ in (1.1). $E$ is equipped with the Euclidean distance $|\cdot|$ and the $l_\infty$-distance $|\cdot|_\infty$. We denote a generic element of $E$ by $x = (u, v)$, $u \in \mathbb{T}_N^d$, $v \in \mathbb{Z}$ and the corresponding closed ball of $|\cdot|_\infty$-radius $r > 0$ centered at $x \in E$ by $B_\infty(x, r)$. Note that $E$ is the image of $\mathbb{Z}^{d+1} = \mathbb{Z}^d \times \mathbb{Z}$ by the mapping $\pi_E \colon \mathbb{Z}^d \times \mathbb{Z} \to E$, $(u, v) \mapsto (\pi_{\mathbb{T}_N^d}(u), v)$, where $\pi_{\mathbb{T}_N^d}$ denotes the canonical projection from $\mathbb{Z}^d$ onto the torus $\mathbb{T}_N^d$. We write $\{e_i\}_{i=1}^{d+1}$ for the canonical basis of $\mathbb{R}^{d+1}$. The projections $\pi_i$, $i = 1, \ldots, d+1$ onto the $d$-dimensional hyperplanes of $E$ are the mappings from $E$ to $(\mathbb{Z}/N\mathbb{Z})^{d-1} \times \mathbb{Z}$ when $i = 1, \ldots, d$, or to $(\mathbb{Z}/N\mathbb{Z})^d$ when $i = d + 1$, defined by omitting the $i$th component of $(u, v) = (u^1, \ldots, u^d, v) \in E$. These projections are not to be confused with the $\mathbb{Z}$-projection $\pi_{\mathbb{Z}}$ from $E$ onto $\mathbb{Z}$,

$$(2.1) \qquad \pi_{\mathbb{Z}}(x) = x \cdot e_{d+1}.$$

For any subset $A \subseteq E$ and $l \geq 1$, we define the $l$-neighborhood of $A$,

$$(2.2) \qquad A^{(l)} = \{x \in E \colon \text{for some } x' \in A, |x - x'|_\infty \leq l\},$$

its $l$-interior,

$$(2.3) \qquad A^{(-l)} = \{x \in A : \text{for all } x' \notin A, |x - x'|_\infty > l\}$$

(so that $A \subseteq B^{(-l)}$ if and only if $A^{(l)} \subseteq B$) and its diameter

$$(2.4) \qquad \text{diam}(A) = \sup\{|x - x'|_\infty : x, x' \in A\}.$$

Given another subset $B \subseteq E$, we define the $B$-relative boundary of $A$,

$$(2.5) \qquad \partial_B(A) = \{x \in B \setminus A : \text{for some } x' \in A, |x - x'| = 1\},$$

and the $B$-relative boundary of $A$ in direction $i \in \{1, \ldots, d+1\}$,

$$(2.6) \qquad \partial_{B,i}(A) = \{x \in B \setminus A : \text{for some } x' \in A, |x - x'| = 1 \text{ and } \pi_i(x) = \pi_i(x')\}.$$

The cube of side-length $l - 1$, $l = 1, \ldots, N$ is defined as

$$(2.7) \qquad C(l) = [0, l-1]^{d+1} \subseteq E$$

(where $[0, l-1] = \{0, \ldots, l-1\}$) and the same cube with base-point $x \in E$ as

$$(2.8) \qquad C_x(l) = x + C(l),$$

where, for $x \in E$ and $A \subset E$, we set $x + A = \{x + x' : x' \in A\} \subseteq E$. For any fixed $N \geq 1$, we now define the probability space

$$(2.9) \qquad \Omega_N = \mathscr{E}([0, \infty), \mathbb{T}_N^d) \times \mathscr{E}([0, \infty), \mathbb{Z}),$$

where, for a set $H$, $\mathscr{E}([0, \infty), H)$ is the space of piecewise constant, right-continuous functions from $[0, \infty)$ to $H$ with infinitely many discontinuities and at most finitely many discontinuities on compact intervals. The canonical processes $(\bar{X}_t)_{t \in [0, \infty)}$, $(\bar{Y}_t)_{t \in [0, \infty)}$ and $(\bar{Z}_t)_{t \in [0, \infty)}$ are defined on $\Omega_N$ by $\bar{X}_t((\omega^{(1)}, \omega^{(2)})) = (\bar{Y}_t, \bar{Z}_t)(\omega^{(1)}, \omega^{(2)}) = (\omega_t^{(1)}, \omega_t^{(2)}) \in E$. These processes generate the canonical filtration $(\bar{\mathcal{F}}_t)_{t \in [0, \infty)}$ on $\Omega_N$ and have the associated shift operators $(\bar{\theta}_t)_{t \in [0, \infty)}$, as well as the jump times $(J_n^{\bar{X}})_{n \geq 0}$, $(J_n^{\bar{Y}})_{n \geq 0}$, $(J_n^{\bar{Z}})_{n \geq 0}$ and counting processes $(N_t^{\bar{X}})_{t \in [0, \infty)}$, $(N_t^{\bar{Y}})_{t \in [0, \infty)}$ and $(N_t^{\bar{Z}})_{t \in [0, \infty)}$, defined for $\bar{X}$ (and analogously for $\bar{Y}$ and $\bar{Z}$) as

$$(2.10) \qquad \begin{aligned} &J_0^{\bar{X}} = 0, \qquad J_1^{\bar{X}} = \inf\{t > 0 : \bar{X}_t \neq \bar{X}_{t-}\} \in (0, \infty), \\ &J_n^{\bar{X}} = J_1^{\bar{X}} \circ \bar{\theta}_{J_{n-1}^{\bar{X}}} + J_{n-1}^{\bar{X}} \qquad \text{for } n \geq 2, \end{aligned}$$

$$(2.11) \qquad N_t^{\bar{X}} = \sup\{n \geq 0 : J_n^{\bar{X}} \leq t\} < \infty, \qquad t \in [0, \infty).$$

The discrete-time processes $(X_n)_{n \geq 0}$, $(Y_n)_{n \geq 0}$ and $(Z_n)_{n \geq 0}$ corresponding to $\bar{X}$, $\bar{Y}$ and $\bar{Z}$ are obtained by restricting time to the integers $n \geq 0$, that is,

$$(2.12) \qquad X_n = \bar{X}_{J_n^{\bar{X}}}, \qquad Y_n = \bar{Y}_{J_n^{\bar{Y}}}, \qquad Z_n = \bar{Z}_{J_n^{\bar{Z}}}, \qquad n \geq 0.$$

Note that, as a consequence one obtains [cf. (2.11)]

$$(2.13) \qquad \bar{X}_t = X_{N_t^{\bar{X}}}, \qquad \bar{Y}_t = Y_{N_t^{\bar{Y}}}, \qquad \bar{Z}_t = Z_{N_t^{\bar{Z}}}, \qquad t \geq 0.$$

For the process $X$, we also define the discrete-time shift operators $(\theta_n)_{n \geq 0}$ and the discrete-time filtration $(\mathcal{F}_n)_{n \geq 0}$ as $\theta_n = \bar{\theta}_{J_n^{\bar{X}}}$, $\mathcal{F}_n = \sigma(X_1, \dots, X_n)$.

We proceed to construct the probability measures $P_x^\Delta$, for $x = (u, v) \in E$ and $0 \leq \Delta < 1$ on $(\Omega_N, (\bar{\mathcal{F}}_t)_{t \in [0, \infty)})$ (and write $E_x^\Delta$ for the corresponding expectations) such that, under $P_x^\Delta$,

$$(2.14) \qquad Y, J^{\bar{Y}}, Z, J^{\bar{Z}} \text{ [cf. (2.10), (2.12)] are independent,}$$

$$(2.15) \qquad Y \text{ is a simple random walk on } \mathbb{T}_N^d \text{ with starting point } u,$$

$Z$ is a random walk on $\mathbb{Z}$ starting at $v$ with transition probability

$$(2.16) \qquad p_Z(v', v' - 1) = \frac{1 - \Delta}{2}, \qquad p_Z(v', v' + 1) = \frac{1 + \Delta}{2}, \qquad v' \in \mathbb{Z}$$

(so $\Delta$ can be interpreted as the drift of the walk in the $\mathbb{Z}$-component),

$$(2.17) \qquad (J_n^{\bar{Y}} - J_{n-1}^{\bar{Y}})_{n \geq 1} \text{ [cf. (2.10)] are i.i.d. Exp(1) variables}$$

[here and throughout this article, $\text{Exp}(\rho)$ denotes the exponential distribution with parameter $\rho > 0$], and

$$(2.18) \qquad (J_n^{\bar{Z}} - J_{n-1}^{\bar{Z}})_{n \geq 1} \text{ are i.i.d. } \text{Exp}\left(\frac{1}{d}\right) \text{ variables.}$$

It follows from this construction that, under $P_x^\Delta$, $X$ is a random walk on $E$ with drift $\Delta$ starting at $x$, that is, a Markov chain on $E$ with initial distribution $\delta_{\{x\}}$ and transition probability specified in (1.2) (in particular, the notation $P_x^\Delta$, $x \in E$, is consistent with its use in the Introduction). Furthermore,

$$(2.19) \qquad (J_n^{\bar{X}} - J_{n-1}^{\bar{X}})_{n \geq 1} \text{ are i.i.d. } \text{Exp}\left(\frac{d+1}{d}\right) \text{ variables}$$

and

$$(2.20) \qquad N^{\bar{X}}, N^{\bar{Y}} \text{ and } N^{\bar{Z}} \text{ [cf. (2.11)] are Poisson processes}$$

on $[0, \infty)$ with respective intensities $\dfrac{d+1}{d}$, $1$ and $\dfrac{1}{d}$.

The disconnection time $T_N^{\text{disc}}$ was defined in (1.3). It will also be useful to consider its continuous-time analog

$$(2.21) \qquad \bar{T}_N^{\text{disc}} = \inf\{t \in [0, \infty) : \bar{X}([0, t]) \text{ disconnects } E\}.$$

Moreover, we will frequently use the following stopping times: The entrance time $H_A^X$ of the set $A \subseteq E$,

$$(2.22) \qquad H_A^X = \inf\{n \geq 0; X_n \in A\},$$

where we write $H_x^X$ if $A = \{x\}$, and the cover time $C_A^X$ of $A \subseteq E$,

$$(2.23) \qquad C_A^X = \inf\{n \geq 0; X([0, n]) \supseteq A\},$$

with obvious modifications such as $H_A^{\tilde{Z}}$ for processes other than $X$ in either discrete or continuous time. For the random walk $X$ and any sets $A \subseteq \bar{A} \subseteq E$, the successive returns $(R_n)_{n \geq 1}$ to $A$ and departures $(D_n)_{n \geq 1}$ from $\bar{A}$ are defined as

$$(2.24) \qquad \begin{aligned} R_1 &= H_A^X, \qquad D_1 = H_{\bar{A}^c}^X \circ \theta_{R_1} + R_1 \quad \text{and} \quad \text{for } n \geq 2, \\ R_n &= R_1 \circ \theta_{D_{n-1}} + D_{n-1}, \qquad D_n = D_1 \circ \theta_{D_{n-1}} + D_{n-1}, \end{aligned}$$

so that $0 \leq R_1 \leq D_1 \leq \cdots \leq R_n \leq D_n \leq \cdots \leq \infty$ and $P_x^\Delta$-a.s. all these inequalities are strict, except possibly the first one. Finally, we also use the Green function of the simple random walk $X$ without drift, killed when exiting $A \subseteq E$, defined as

$$(2.25) \qquad g^A(x, x') = E_x^0\left[\sum_{n=0}^{\infty} \mathbf{1}\{X_n = x', n < H_{A^c}^X\}\right], \qquad x, x' \in E.$$

We conclude this section with the Girsanov-type estimate comparing $P_x^\Delta$ and $P_x^0$.

LEMMA 2.1 $[d \geq 1, \ N \geq 1, \ \Delta \in (0, 1), \ x \in E]$. *Consider any $(\mathcal{F}_n)_{n \geq 0}$-stopping time $T$ and any $\mathcal{F}_T$-measurable event $A$ such that, for some $b, b' \in \mathbb{R} \cup \{-\infty, \infty\}$,*

$$(2.26) \qquad T < \infty \quad \text{and} \quad b \leq \pi_{\mathbb{Z}}(X_T - x) \leq b', \qquad P_x^0\text{-a.s. on } A.$$

*Then*

$$(2.27) \qquad (1 - \Delta)^{b_-}(1 + \Delta)^{b_+} E_x^0[A, (1 - \Delta^2)^{[T/2]}] \leq P_x^\Delta[A]$$

*and*

$$(2.28) \qquad P_x^\Delta[A] \leq (1 - \Delta)^{b'_-}(1 + \Delta)^{b'_+} P_x^0[A],$$

*where we set $(1 - \Delta)^\infty = 0$ and $(1 + \Delta)^\infty = \infty$.*

PROOF. For any $\mathcal{F}_n$-measurable event $A_n$, it follows directly from the definition of the transition probabilities of the walk $X$ [cf. (1.2)] that

$$(2.29) \qquad P_x^\Delta[A_n] = E_x^0\left[A_n, \prod_{i=1}^{n}(1 + \Delta\pi_{\mathbb{Z}}(X_i - X_{i-1}))\right].$$

For any $(\mathcal{F}_n)_{n \geq 0}$-stopping time $T$ satisfying (2.26), we apply (2.29) with the $\mathcal{F}_n$-measurable event $A_n = A \cap \{T = n\}$ for $n \geq 0$ and deduce, via monotone convergence,

$$
\begin{aligned}
(2.30) \qquad P_x^\Delta[A] &= \sum_{n \geq 0} P_x^\Delta[A_n] = \sum_{n \geq 0} E_x^0 \left[ A_n, \prod_{i=1}^{T} (1 + \Delta \pi_{\mathbb{Z}}(X_i - X_{i-1})) \right] \\
&= E_x^0 \left[ A, \prod_{i=1}^{T} (1 + \Delta \pi_{\mathbb{Z}}(X_i - X_{i-1})) \right].
\end{aligned}
$$

To complete the proof, we bound the product inside the expectation on the right-hand side of (2.30) from above and from below. The contribution of the product is a factor of $1 + \Delta$ for every displacement of $X$ into the positive $\mathbb{Z}$-direction up to time $T$ and a factor of $1 - \Delta$ for every displacement into the negative $\mathbb{Z}$-direction during the same time. We now group together the factors in the product as pairs of the form $(1 + \Delta)(1 - \Delta) = 1 - \Delta^2$ for as many factors as possible (i.e., until all remaining factors are of the form $1 + \Delta$ or all remaining factors are of the form $1 - \Delta$). By (2.26), the contribution of these remaining factors is bounded from below by $(1 - \Delta)^{b_-}(1 + \Delta)^{b_+}$ and from above by $(1 - \Delta)^{b'_-}(1 + \Delta)^{b'_+}$. For (2.28), we note that $1 - \Delta^2 < 1$ and bound the contribution made by the pairs from above by 1. For (2.27), we note that the number of pairs contributed can be at most $[\frac{T}{2}]$. This completes the proof of the lemma. $\quad \square$

## 3. Upper bounds.

This section is devoted to upper bounds on $T_N^{\mathrm{disc}}$. We will prove the following theorem, which is more than sufficient to yield the upper bounds in Theorem 1.1:

THEOREM 3.1 $(d \geq 2,\ \alpha > 0,\ \varepsilon > 0)$. *For some constant $c_0 > 0$,*

$$
(3.1) \qquad \text{for } \alpha > 1, \qquad\qquad P_0^{N^{-d\alpha}}[T_N^{\mathrm{disc}} \leq N^{2d}(\log N)^{4+\varepsilon}] \overset{N \to \infty}{\longrightarrow} 1,
$$

$$
(3.2) \qquad \text{for } \alpha \leq 1, \qquad P_0^{N^{-d\alpha}}[T_N^{\mathrm{disc}} \leq \exp\{c_0 N^{d(1-\alpha)}(\log N)^2\}] \overset{N \to \infty}{\longrightarrow} 1.
$$

In order to show Theorem 3.1, it suffices to show the corresponding result in continuous time, which is [cf. (2.21)]:

THEOREM 3.2 $(d \geq 2,\ \alpha > 0,\ \varepsilon > 0)$. *For some constant $c_0 > 0$,*

$$
(3.3) \qquad \text{for } \alpha > 1, \qquad\qquad P_0^{N^{-d\alpha}}[\bar{T}_N^{\mathrm{disc}} \leq N^{2d}(\log N)^{4+\varepsilon}] \overset{N \to \infty}{\longrightarrow} 1,
$$

$$
(3.4) \qquad \text{for } \alpha \leq 1, \qquad P_0^{N^{-d\alpha}}[\bar{T}_N^{\mathrm{disc}} \leq \exp\{c_0 N^{d(1-\alpha)}(\log N)^2\}] \overset{N \to \infty}{\longrightarrow} 1.
$$

PROOF THAT THEOREM 3.2 IMPLIES THEOREM 3.1. For this proof as well as for future reference, we note that, for any Poisson process $(N_t^{(\rho)})_{t \in [0,\infty)}$ of parameter $\rho > 0$, one has by the exponential Chebyshev inequality, for any $t \geq 0$,

$$P[N_t^{(\rho)} \geq e\rho t] \leq e^{-e\rho t} E[e^{N_t^{(\rho)}}] = e^{-e\rho t + \rho t(e-1)} = e^{-\rho t},$$

as well as

$$P[N_t^{(\rho)} \leq e^{-1}\rho t] \leq e^{e^{-1}\rho t} E[e^{-N_t^{(\rho)}}] = e^{e^{-1}\rho t + \rho t(e^{-1}-1)} = e^{-\rho t(1-2e^{-1})},$$

hence

$$(3.5) \qquad P[N_t^{(\rho)} \notin (e^{-1}\rho t, e\rho t]] \leq 2e^{-c\rho t}.$$

Let us now assume that Theorem 3.2 is true. By definition of $X$ and $\bar{X}$ [cf. (2.13)], one has, for any $s, t \geq 0$, on the event $\{N_t^{\bar{X}} \leq s\}$,

$$\{T_N^{\mathrm{disc}} > s\} = \{X([0,[s]]) \text{ does not disconnect } E\}$$

$$\subseteq \{X([0, N_t^{\bar{X}}]) \text{ does not disconnect } E\} \overset{((2.13),(2.21))}{=} \{\bar{T}_N^{\mathrm{disc}} > t\}.$$

Using this last observation with $s = e\frac{d+1}{d}t$, we deduce

$$(3.6) \qquad \begin{aligned}
P_0^{N^{-d\alpha}}\left[T_N^{\mathrm{disc}} > e\frac{d+1}{d}t\right] &\leq P_0^{N^{-d\alpha}}\left[T_N^{\mathrm{disc}} > e\frac{d+1}{d}t, N_t^{\bar{X}} \leq e\frac{d+1}{d}t\right] \\
&\quad + P_0^{N^{-d\alpha}}\left[T_N^{\mathrm{disc}} > e\frac{d+1}{d}t, N_t^{\bar{X}} > e\frac{d+1}{d}t\right] \\
&\leq P_0^{N^{-d\alpha}}[\bar{T}_N^{\mathrm{disc}} > t] + P_0^{N^{-d\alpha}}\left[N_t^{\bar{X}} > e\frac{d+1}{d}t\right].
\end{aligned}$$

We now fix any $\alpha > 1$ and $\varepsilon > 0$. The last inequality with $t_N = N^{2d}(\log N)^{4+\varepsilon/2}$ yields, for $N \geq c(\varepsilon)$ (we refer to the end of the [Introduction](#) for our convention concerning constants),

$$\begin{aligned}
&P_0^{N^{-d\alpha}}[T_N^{\mathrm{disc}} \geq N^{2d}(\log N)^{4+\varepsilon}] \\
&\qquad \leq P_0^{N^{-d\alpha}}\left[T_N^{\mathrm{disc}} > e\frac{d+1}{d}t_N\right] \\
&\qquad \overset{(3.6)}{\leq} P_0^{N^{-d\alpha}}[\bar{T}_N^{\mathrm{disc}} > t_N] + P_0^{N^{-d\alpha}}\left[N_{t_N}^{\bar{X}} > e\frac{d+1}{d}t_N\right].
\end{aligned}$$

The first of the two terms on the right-hand side tends to 0 as $N \to \infty$, by (3.3), while the second term is bounded from above by $2\exp\{-cN^{2d}(\log N)^{4+\varepsilon/2}\}$ by (3.5). We have thus deduced (3.1).

For $\alpha \leq 1$, we proceed in the same way: Applied with $t'_N = \exp\{c_0 N^{d(1-\alpha)} \times (\log N)^2\}$, (3.6) and (3.5) yield, for $N \geq c(c_0)$,

$$P_0^{N^{-d\alpha}}[T_N^{\mathrm{disc}} \geq \exp\{2c_0 N^{d(1-\alpha)}(\log N)^2\}]$$

$$\leq P_0^{N^{-d\alpha}}\left[T_N^{\mathrm{disc}} > e\frac{d+1}{d}t'_N\right]$$

$$\leq P_0^{N^{-d\alpha}}[\bar{T}_N^{\mathrm{disc}} > t'_N] + \exp\{-ce^{c_0 N^{(1-\alpha)}(\log N)^2}\},$$

so that (3.2) follows from (3.4). $\quad\square$

PROOF OF THEOREM 3.2. Following the idea outlined in the Introduction, we define the process $V$, whose purpose is to record visits of $X$ to $\mathbb{T}_N^d \times \{0\}$. To this end, we introduce the stopping times $(\bar{S}_n)_{n\geq 0}$ by setting [cf. (2.10), (2.22)]

(3.7)
$$\bar{S}_0 = 0, \qquad \bar{S}_1 = H_0^{\bar{Z}} \circ \bar{\theta}_{J_1^{\bar{Y}}} + J_1^{\bar{Y}} \leq \infty, \qquad \text{and for } n \geq 2,$$

$$\bar{S}_n = \begin{cases} \bar{S}_1 \circ \bar{\theta}_{\bar{S}_{n-1}} + \bar{S}_{n-1}, & \text{on } \{\bar{S}_{n-1} < \infty\}, \\ \infty, & \text{on } \{\bar{S}_{n-1} = \infty\}, \end{cases}$$

and on the event $\{\bar{S}_k < \infty\}$, we define

(3.8)
$$V_n = \bar{Y}_{\bar{S}_n}, \qquad n = 0, \dots, k.$$

Note that, as soon as $V$ has visited all points of $\mathbb{T}_N^d$, $\bar{X}$ has visited all points of $\mathbb{T}_N^d \times \{0\}$, and has therefore disconnected $E$. Hence, we are interested in an upper bound on the cover time $C_{\mathbb{T}_N^d}^V$ [cf. (2.23)]. This desired upper bound will result from the following estimate on cover times for symmetric Markov chains. Following Aldous and Fill ([1], Chapter 7, page 2), we call a Markov chain $(W_n)_{n\geq 0}$ on the finite state-space $G$ with transition probabilities $p_W(g, g')$, $g, g' \in G$ *symmetric*, if for any states $g_0, g_1 \in G$, there exists a bijection $\gamma : G \to G$ satisfying $\gamma(g_0) = g_1$ and $p_W(g, g') = p_W(\gamma(g), \gamma(g'))$ for all $g, g' \in G$.

LEMMA 3.3. *Given a symmetric, irreducible and reversible Markov chain $(W_n)_{n\geq 0}$ on the finite state-space $G$ whose transition matrix $(p_W(g, g'))_{g, g' \in G}$ has eigenvalues $1 = \lambda_1(W) > \lambda_2(W) \geq \cdots \geq \lambda_{|G|}(W) \geq -1$, one has*

(3.9)
$$P_g[C_G^W \geq n] \leq |G| \exp\left\{-\left[\frac{n}{4eu(W)}\right]\right\} \qquad \text{for any } g \in G, \ n \geq 1,$$

*where $P_g$ is the canonical probability on $G^{\mathbb{N}}$ governing $W$ with $W_0 = g$ and*

(3.10)
$$u(W) = \sum_{m=2}^{|G|} \frac{1}{1 - \lambda_m(W)}.$$

Proof. We assume that $n \geq 4eu(W)$, for otherwise there is nothing to prove. The following estimate on the maximum hitting time [cf. (2.22)] is a consequence of the so-called eigentime identity (see [1], Lemma 15 and Proposition 13 in Chapter 3, and note that $E_g[H_{g'}^W] = E_{g'}[H_g^W]$ by our assumptions on symmetry, irreducibility and reversibility; cf. [1], Chapter 3, Lemma 1):

$$(3.11) \qquad \max_{g,g' \in G} E_g[H_{g'}^W] \leq 2 \sum_{m=2}^{|G|} \frac{1}{1 - \lambda_m(W)} \overset{(3.10)}{=} 2u(W).$$

Choosing any $1 \leq s \leq n$, we deduce the following tail estimate on $C_G^W$ with a standard application of the simple Markov property at the times $([s] - 1)[\frac{n}{s}], \ldots, 2[\frac{n}{s}], [\frac{n}{s}]$:

$$\begin{aligned}
P_g[&C_G^W \geq n] \\
&= P_g[\text{for some } g' \in G : H_{g'}^W \geq n] \\
&\leq |G| \max_{g,g' \in G} P_g[H_{g'}^W \geq n] \overset{(\text{Markov})}{\leq} |G| \left( \max_{g,g' \in G} P_g \left[ H_{g'}^W \geq \left[ \frac{n}{s} \right] \right] \right)^{[s]} \\
&\overset{(\text{Chebyshev, } (3.11))}{\leq} |G| \left( \left[ \frac{n}{s} \right]^{-1} 2u(W) \right)^{[s]} \overset{(n/s \leq 2[n/s])}{\leq} |G| \left( \frac{4su(W)}{n} \right)^{[s]}.
\end{aligned}$$

With $1 \leq s = \frac{n}{4eu(W)} \leq n$, this yields (3.9).  □

In what follows, we require the following alternative expression for the distribution of the stopping times $(\bar{S}_n)_{n \geq 0}$ [cf. 3.7]:

Lemma 3.4 ($d \geq 1$, $n \geq 1$, $\Delta \in [0,1)$). *The following equality in distribution holds under $P_0^\Delta$:*

$$(3.12) \qquad \bar{S}_n \overset{(\text{dist.})}{=} \sigma_1 + \cdots + \sigma_n + H_{-\hat{Z}_{\sigma_1}^{(1)}}^{\bar{Z}^{(1)}} + \cdots + H_{-\hat{Z}_{\sigma_n}^{(n)}}^{\bar{Z}^{(n)}},$$

*where the random variables $\{\sigma_i\}_{i \geq 1}$ and the processes $\{\bar{Z}^{(i)}, \hat{Z}^{(i)}\}_{i \geq 1}$ are independent, $\{\bar{Z}^{(i)}, \hat{Z}^{(i)}\}_{i \geq 1}$ are i.i.d. copies of the random walk $\bar{Z}$ and the $\sigma_i$ are exponentially distributed with parameter 1.*

Proof. It suffices to prove that, for any $n \geq 1$,

$$(3.13) \quad (\bar{S}_1, \bar{S}_2 - \bar{S}_1, \ldots, \bar{S}_n - \bar{S}_{n-1}) \overset{(\text{dist.})}{=} (\sigma_1 + H_{-\hat{Z}_{\sigma_1}^{(1)}}^{\bar{Z}^{(1)}}, \ldots, \sigma_n + H_{-\hat{Z}_{\sigma_n}^{(n)}}^{\bar{Z}^{(n)}}),$$

for then we obtain

$$\bar{S}_n = \sum_{i=1}^n (\bar{S}_i - \bar{S}_{i-1}) \overset{(\text{dist.})}{=} \sum_{i=1}^n (\sigma_i + H_{-\hat{Z}_{\sigma_i}^{(n)}}^{\bar{Z}^{(i)}}),$$

as required. For the purpose of showing (3.13), we fix any $t_1, \ldots, t_n \geq 0$ and find, with the strong Markov property:

$$(3.14) \quad \begin{aligned} P_0^\Delta \left[ \bigcap_{i=1}^n \{ \bar{S}_i - \bar{S}_{i-1} \leq t_i \} \right] &= P_0^\Delta \left[ \bigcap_{i=1}^{n-1} \{ \bar{S}_i - \bar{S}_{i-1} \leq t_i \} \cap \theta_{\bar{S}_{n-1}}^{-1} \{ \bar{S}_1 \leq t_n \} \right] \\ &= E_0^\Delta \left[ \bigcap_{i=1}^{n-1} \{ \bar{S}_i - \bar{S}_{i-1} \leq t_i \}, P_{\bar{X}_{\bar{S}_{n-1}}}^\Delta [\bar{S}_1 \leq t_n] \right]. \end{aligned}$$

Thanks to translation invariance in the $\mathbb{T}_N^d$-direction, the distribution of $\bar{S}_n$, $n \geq 0$, does not depend on the $\mathbb{T}_N^d$-coordinate of the starting point. In particular, one has

$$(3.15) \qquad P_{(u,0)}^\Delta [\bar{S}_n \leq \cdot] = P_0^\Delta [\bar{S}_n \leq \cdot] \qquad \text{for any } u \in \mathbb{T}_N^d, n \geq 0.$$

Therefore (3.14) simplifies to

$$(3.16) \quad \begin{aligned} &P_0^\Delta \left[ \bigcap_{i=1}^n \{ \bar{S}_i - \bar{S}_{i-1} \leq t_i \} \right] \\ &= P_0^\Delta \left[ \bigcap_{i=1}^{n-1} \{ \bar{S}_i - \bar{S}_{i-1} \leq t_i \} \right] P_0^\Delta [\bar{S}_1 \leq t_n] \\ &\overset{\text{(induction)}}{=} \prod_{i=1}^n P_0^\Delta [\bar{S}_1 \leq t_i]. \end{aligned}$$

However, $J_1^{\bar{Y}}$ is exponentially distributed with parameter 1 [cf. (2.17)] and independent of $\bar{Z}$ [cf. (2.14)]. We hence obtain by Fubini's theorem that, for any $t \geq 0$:

$$(3.17) \quad \begin{aligned} P_0^\Delta [\bar{S}_1 \leq t] &\overset{(3.7)}{=} P_0^\Delta [H_0^{\bar{Z}} \circ \bar{\theta}_{J_1^{\bar{Y}}} + J_1^{\bar{Y}} \leq t] \\ &\overset{\text{(Fub.)}}{=} \int_0^t P_0^\Delta [H_0^{\bar{Z}} \circ \bar{\theta}_s + s \leq t] e^{-s} \, ds. \end{aligned}$$

Applying the simple Markov property at time $s$ to the probability inside this last integral, one finds

$$\begin{aligned} P_0^\Delta [H_0^{\bar{Z}} \circ \bar{\theta}_s + s \leq t] &= E_0^\Delta [P_{(0, \bar{Z}_s)}^\Delta [H_0^{\bar{Z}} \leq t - s]] \\ &\overset{(\bar{Z} \overset{\text{(dist.)}}{=} \hat{Z}^{(1)})}{=} E_0^\Delta [P_{(0, \hat{Z}_s^{(1)})}^\Delta [H_0^{\bar{Z}} \leq t - s]] \\ &\overset{\text{(transl. inv.)}}{=} P_0^\Delta [H_{-\hat{Z}_s^{(1)}}^{\bar{Z}} \leq t - s]. \end{aligned}$$

Inserting this last expression into (3.17) and applying again Fubini's theorem, we obtain

$$P_0^\Delta[\bar{S}_1 \le t] = P_0^\Delta[\sigma_1 + H_{-\hat{Z}_{\sigma_1}^{(1)}}^{\bar{Z}} \le t].$$

By this observation and independence of $\{\sigma_i, \bar{Z}^{(i)}, \hat{Z}^{(i)}\}_{i \ge 1}$, (3.16) becomes

$$P_0^\Delta\left[\bigcap_{i=1}^n \{\bar{S}_i - \bar{S}_{i-1} \le t_i\}\right] = P_0^\Delta\left[\bigcap_{i=1}^n \{\sigma_i + H_{-\hat{Z}_{\sigma_i}^{(i)}}^{\bar{Z}^{(i)}} \le t_i\}\right],$$

which shows (3.13) and hence completes the proof of Lemma 3.4.  □

The next step toward the application of Lemma 3.3 is to show that $(V_n)_{n=1}^k$ [cf. (3.8)] satisfies the hypotheses imposed on $W$, provided we take the event $\{\bar{S}_k < \infty\}$ as probability space, equipped with the probability measure $P_{(u,0)}^\Delta(\cdot | \bar{S}_k < \infty)$, $u \in \mathbb{T}_N^d$.

LEMMA 3.5 ($d \ge 1$, $k \ge 1$, $\Delta \in [0,1)$, $u \in \mathbb{T}_N^d$).   On the probability space $(\{\bar{S}_k < \infty\}, P_{(u,0)}^\Delta[\cdot | \bar{S}_k < \infty])$ and the finite time interval $n = 0, \ldots, k$, $(V_n)_{n=0}^k$ is a symmetric, irreducible and reversible Markov chain on $\mathbb{T}_N^d$ starting at $u$ with transition probability

$$(3.18) \qquad p_V(u,u') = P_{(u,0)}^\Delta[\bar{Y}_{\bar{S}_1} = u' | \bar{S}_1 < \infty], \qquad u,u' \in \mathbb{T}_N^d.$$

PROOF.   By construction $Y$, $J^{\bar{Y}}$, $\bar{Z}$ are independent [cf. (2.14)]. Since $\bar{S}_1$ and $N_{\bar{S}_1}^{\bar{Y}}$ are both $\sigma(J^{\bar{Y}}, \bar{Z})$-measurable [cf. (2.11), (3.7)], it follows that $Y$ and $(\bar{S}_1, N_{\bar{S}_1}^{\bar{Y}})$ are independent as well. Hence, one can rewrite the expression for $p_V(u,u')$ in (3.18) using Fubini's theorem:

$$
\begin{aligned}
p_V(u,u') &\overset{(3.15)}{=} \frac{1}{P_0^\Delta[\bar{S}_1 < \infty]} P_{(u,0)}^\Delta[Y_{N_{\bar{S}_1}^{\bar{Y}}} = u', \bar{S}_1 < \infty] \\
(3.19) \qquad &\overset{\text{(Fubini)}}{=} \frac{1}{P_0^\Delta[\bar{S}_1 < \infty]} E_{(u,0)}^\Delta[P_{(u,0)}^\Delta[Y_n = u']|_{n=N_{\bar{S}_1}^{\bar{Y}}}, \bar{S}_1 < \infty] \\
&= \frac{1}{P_0^\Delta[\bar{S}_1 < \infty]} E_0^\Delta[P_{(u,0)}^\Delta[Y_n = u']|_{n=N_{\bar{S}_1}^{\bar{Y}}}, \bar{S}_1 < \infty],
\end{aligned}
$$

where in the last line we have used that the expression inside the expectation is a function of $N_{\bar{S}_1}^{\bar{Y}}$ and $\bar{S}_1$ and therefore does not depend on the $\mathbb{T}_N^d$-coordinate of the starting point. From (3.19), it follows that the transition probabilities $p_V(\cdot, \cdot)$ define an irreducible, symmetric (as defined above Lemma 3.3) and reversible process. Indeed, for any $u, u' \in \mathbb{T}_N^d$ such that $P_{(u,0)}^\Delta[Y_1 = u'] > 0$, (3.19) and $P_0^\Delta[N_{\bar{S}_1}^{\bar{Y}} = 1, \bar{S}_1 < \infty] \ge P_0^\Delta[X_1 \in \mathbb{T}_N^d \times \{0\}] >$ 0 imply that $p_V(u, u') > 0$, so that irreducibility follows from irreducibility of the simple random walk $Y$. Similarly, (3.19) shows that symmetry follows from symmetry of $Y$, which holds by translation invariance. Finally, reversibility follows by exchanging $u$ and $u'$ in the last line of (3.19), which one can do by reversibility of $Y$. It thus remains to be shown that $p_V(\cdot, \cdot)$ are in fact the correct transition probabilities for $V$, that is, that for any $u, u_1, \ldots, u_n \in \mathbb{T}_N^d$, $1 \le n \le k$, and

$$(3.20) \qquad\qquad A = \{V_0 = u, \ldots, V_{n-1} = u_{n-1}\},$$

one has

$$(3.21) \qquad P_{(u,0)}^\Delta[V_n = u_n, A | \bar{S}_k < \infty] = p_V(u_{n-1}, u_n) P_{(u,0)}^\Delta[A | \bar{S}_k < \infty].$$

Using the strong Markov property at time $\bar{S}_n$, one has

$$P_{(u,0)}^\Delta[V_n = u_n, A | \bar{S}_k < \infty]$$

$$
\begin{aligned}
(3.22) \quad &\overset{\text{(Markov)}}{=} \frac{1}{P_0^\Delta[\bar{S}_k < \infty]} \\
&\qquad \times E_{(u,0)}^\Delta[\bar{Y}_{\bar{S}_n} = u_n, A, \bar{S}_n < \infty, P_{(\bar{Y}_{\bar{S}_n}, 0)}^\Delta[\bar{S}_{k-n} < \infty]] \\
&\overset{(3.15)}{=} \frac{P_0^\Delta[\bar{S}_{k-n} < \infty]}{P_0^\Delta[\bar{S}_k < \infty]} P_{(u,0)}^\Delta[\bar{Y}_{\bar{S}_n} = u_n, A, \bar{S}_n < \infty].
\end{aligned}
$$

Applying the strong Markov property at time $\bar{S}_{n-1}$ to the last probability in this expression, we infer that

$$P_{(u,0)}^\Delta[\bar{Y}_{\bar{S}_n} = u_n, A, \bar{S}_n < \infty]$$

$$
\begin{aligned}
&\overset{\text{(Markov)}}{=} E_{(u,0)}^\Delta[A, \bar{S}_{n-1} < \infty, P_{(\bar{Y}_{\bar{S}_{n-1}}, 0)}^\Delta[\bar{Y}_{\bar{S}_1} = u_n, \bar{S}_1 < \infty]] \\
&\overset{((3.18), (3.20))}{=} p_V(u_{n-1}, u_n) E_{(u,0)}^\Delta[A, \bar{S}_{n-1} < \infty, P_{(u_{n-1},0)}^\Delta[\bar{S}_1 < \infty]] \\
&\overset{(3.20)}{=} p_V(u_{n-1}, u_n) E_{(u,0)}^\Delta[A, \bar{S}_{n-1} < \infty, P_{\bar{X}_{\bar{S}_{n-1}}}^\Delta[\bar{S}_1 < \infty]] \\
&\overset{\text{(Markov)}}{=} p_V(u_{n-1}, u_n) P_{(u,0)}^\Delta[A, \bar{S}_n < \infty].
\end{aligned}
$$

Substituting this last expression into (3.22), and noting that (once more by the strong Markov property)

$$
\begin{aligned}
P_0^\Delta[\bar{S}_{k-n} < \infty] P_{(u,0)}^\Delta[A, \bar{S}_n < \infty] &\overset{(3.15)}{=} E_{(u,0)}^\Delta[A, \bar{S}_n < \infty, P_{\bar{X}_{\bar{S}_n}}^\Delta[\bar{S}_{k-n} < \infty]] \\
&\overset{\text{(Markov)}}{=} P_{(u,0)}^\Delta[A, \bar{S}_k < \infty],
\end{aligned}
$$

we obtain (3.21) and finish the proof of Lemma 3.5.   $\square$

With the notation of Lemma 3.3, we recall that $\lambda_m(V)$ and $\lambda_m(Y)$, $m = 1, \ldots, N^d$ stand for the decreasingly ordered eigenvalues of the transition matrices $(p_V(u, u'))_{u, u' \in \mathbb{T}_N^d}$ and $(p_Y(u, u'))_{u, u' \in \mathbb{T}_N^d}$ of $V$ and $Y$, respectively. The following proposition shows how these two sets of eigenvalues are related.

PROPOSITION 3.6 ($d \geq 1$, $\Delta \in [0, 1)$).

$$(3.23) \qquad \lambda_m(V) = E_0^\Delta[\lambda_m(Y)^{N_{\bar{S}_1}^{\bar{Y}}} | \bar{S}_1 < \infty], \qquad 1 \leq m \leq N^d.$$

PROOF. From (3.19), we know that, for $u$, $u' \in \mathbb{T}_N^d$,

$$p_V(u, u') = E_0^\Delta[p_Y^{N_{\bar{S}_1}^{\bar{Y}}}(u, u') | \bar{S}_1 < \infty].$$

For any eigenvalue/eigenvector pair $(\lambda_m(Y), v_m)$, we infer that

$$(p_V(u, u'))_{u, u'} v_m = E_0^\Delta[(p_Y(u, u'))_{u, u'}^{N_{\bar{S}_1}^{\bar{Y}}} v_m | \bar{S}_1 < \infty]$$
$$= E_0^\Delta[\lambda_m(Y)^{N_{\bar{S}_1}^{\bar{Y}}} v_m | \bar{S}_1 < \infty] = E_0^\Delta[\lambda_m(Y)^{N_{\bar{S}_1}^{\bar{Y}}} | \bar{S}_1 < \infty] v_m.$$

Hence, $(p_V(u, u'))_{u, u' \in \mathbb{T}_N^d}$ has the same eigenvectors as $(p_Y(u, u'))_{u, u' \in \mathbb{T}_N^d}$ and the corresponding eigenvalues are indeed given by (3.23). □

We can thus relate the quantity $u(V)$ to $u(Y)$ [cf. (3.10)], which is well known from Aldous and Fill [1]:

PROPOSITION 3.7 ($d \geq 2$, $N \geq 1$).

$$(3.24) \qquad u(Y) \leq cN^2 \log N \qquad (d = 2),$$

$$(3.25) \qquad u(Y) \leq cN^d \qquad (d \geq 3).$$

(We refer to the end of the Introduction for our convention concerning constants.)

PROOF. The proof is contained in [1]: By the eigentime identity from Chapter 3, Proposition 13, $u(Y)$ is equal to the average hitting time (cf. Chapter 4, page 1, for the definition), for which the estimates hold by Proposition 8 in Chapter 13. □

As a consequence, we now obtain our desired estimate on $C_{\mathbb{T}_N^d}^V$ by an application of Lemma 3.3:

LEMMA 3.8 ($d \geq 2$, $N \geq 2$, $u \in \mathbb{T}_N^d$). *For any $k \geq [c_1 N^d (\log N)^2]$, one has*

$$(3.26) \qquad \sup_{\Delta \in [0,1)} P_{(u,0)}^{\Delta} [C_{\mathbb{T}_N^d}^V \geq [c_1 N^d (\log N)^2] | \bar{S}_k < \infty] \leq \frac{1}{N^{10}}.$$

PROOF. We fix any $\Delta \in [0,1)$ and consider the canonical Markov chain $(W_n)_{n \geq 0}$, with state-space $\mathbb{T}_N^d$, starting point $u$ and with the same transition probability as $(V_n)_{n=0}^k$ under $P_{u,0}^{\Delta}[\cdot | \bar{S}_k < \infty]$, that is, $p_W(\cdot, \cdot) = p_V(\cdot, \cdot)$. By Lemma 3.5, $(W_n)_{n \geq 0}$ then satisfies the assumptions of Lemma 3.3. Moreover, $(W_n)_{n=0}^k$ has the same distribution as $(V_n)_{n=0}^k$ under $P_{u,0}^{\Delta}[\cdot | \bar{S}_k < \infty]$. With the help of Lemma 3.3, we see that, for $k \geq [cN^d (\log N)^2]$,

$$P_{(u,0)}^{\Delta} [C_{\mathbb{T}_N^d}^V \geq [cN^d (\log N)^2] | \bar{S}_k < \infty]$$

$$(3.27) \qquad = P_u [C_{\mathbb{T}_N^d}^W \geq [cN^d (\log N)^2]]$$

$$\overset{(3.9)}{\leq} N^d \exp \left\{ - \left[ \frac{[cN^d (\log N)^2]}{4eu(W)} \right] \right\}.$$

Since $V$ and $W$ have the same transition probability, we have $u(W) = u(V)$, so once we show that

$$(3.28) \qquad u(V) = \sum_{m=2}^{N^d} \frac{1}{1 - \lambda_m(V)} \leq cN^d + u(Y),$$

the proof of (3.26) will be complete with (3.24), (3.25), (3.27) by choosing $c = c_1$ a large enough constant and noting that the right-hand side of (3.27) does not depend on $\Delta$. We use the expression for $\lambda_m(V)$ of (3.23) and distinguish the two cases $0 < \lambda_m(Y) < 1$ and $-1 \leq \lambda_m(Y) \leq 0$. If $0 < \lambda_m(Y) < 1$, then $\lambda_m(V) \leq \lambda_m(Y)$, because $N_{\bar{S}_1}^{\bar{Y}} \geq 1$ by definition of $\bar{S}_1$ [cf. (3.7)], and hence

$$(3.29) \qquad \frac{1}{1 - \lambda_m(V)} \leq \frac{1}{1 - \lambda_m(Y)}.$$

If, on the other hand, $-1 \leq \lambda_m(Y) \leq 0$, then $\lambda_m(Y)^n$ is nonnegative only for even $n \geq 1$ and not larger than 1 for all $n \geq 1$, so in particular

$$(3.30) \quad \lambda_m(V) \leq P_0^{\Delta} [N_{\bar{S}_1}^{\bar{Y}} \geq 2 | \bar{S}_1 < \infty] \overset{(N_{\bar{S}_1}^{\bar{Y}} \geq 1)}{=} 1 - P_0^{\Delta} [N_{\bar{S}_1}^{\bar{Y}} = 1 | \bar{S}_1 < \infty].$$

Since $\{X_1 \in \mathbb{T}_N^d \times \{0\}\} \subseteq \{\bar{S}_1 = J_1^{\bar{Y}}\} \subseteq \{N_{\bar{S}_1}^{\bar{Y}} = 1\}$ [cf. (3.7)], we deduce from (3.30) that

$$\lambda_m(V) \leq 1 - P_0^{\Delta} [X_1 \in \mathbb{T}_N^d \times \{0\} | \bar{S}_1 < \infty] = 1 - \frac{P_0^{\Delta} [X_1 \in \mathbb{T}_N^d \times \{0\}]}{P_0^{\Delta} [\bar{S}_1 < \infty]}$$

$$\leq 1 - P_0^{\Delta} [X_1 \in \mathbb{T}_N^d \times \{0\}] = 1 - \frac{d}{d+1},$$

and hence

$$(3.31) \qquad \frac{1}{1 - \lambda_m(V)} \leq \frac{d+1}{d}.$$

The estimates (3.29) and (3.31) together yield (3.28), so the proof of Lemma 3.8 is complete. $\square$

In view of (3.26), we still need an upper bound on the amount of time it takes for the corresponding $[c_1 N^d (\log N)^2]$ returns to occur. For simplicity of notation, we set

$$(3.32) \qquad a_N = N^d (\log N)^2,$$

and we treat the cases $\alpha > 1$ and $\alpha \leq 1$ in Theorem 3.2 separately.

*Case $\alpha > 1$.* We observe that

$$(3.33) \qquad P_0^{N-d\alpha}[\bar{T}_N^{\mathrm{disc}} \geq a_N^2 (\log N)^\varepsilon] \leq P_0^{N-d\alpha}[\bar{T}_N^{\mathrm{disc}} \geq \bar{S}_{[c_1 a_N]}]$$
$$+ P_0^{N-d\alpha}[\bar{S}_{[c_1 a_N]} > a_N^2 (\log N)^\varepsilon].$$

By Lemma 3.8 one has

$$P_0^{N-d\alpha}[\bar{T}_N^{\mathrm{disc}} \geq \bar{S}_{[c_1 a_N]}] \leq P_0^{N-d\alpha}[C_{\mathbb{T}_N^d \times \{0\}}^{\tilde{X}} \geq \bar{S}_{[c_1 a_N]}]$$
$$= P_0^{N-d\alpha}[C_{\mathbb{T}_N^d}^V \geq [c_1 a_N], \bar{S}_{[c_1 a_N]} < \infty] \overset{(3.26)}{\longrightarrow} 0.$$

In view of (3.33), the proof of (3.3) will thus be complete once it is shown that

$$(3.34) \qquad P_0^{N-d\alpha}[\bar{S}_{[c_1 a_N]} \leq a_N^2 (\log N)^\varepsilon] \overset{N \to \infty}{\longrightarrow} 1.$$

With (3.12), this will follow from (we refer to the statement of Lemma 3.4 for the notation)

$$(3.35) \qquad P_0^{N-d\alpha}\left[\sum_{i=1}^{[c_1 a_N]} (\sigma_i + H_{-\hat{Z}_{\sigma_i}^{(i)}}^{\bar{Z}^{(i)}}) \leq a_N^2 (\log N)^\varepsilon\right] \overset{N \to \infty}{\longrightarrow} 1.$$

Let us define the event $A(c_1 a_N)$ by

$$(3.36) \qquad A(c_1 a_N) = \{\sigma_1 + \cdots + \sigma_{[c_1 a_N]} \leq [2 c_1 a_N],$$
$$|\hat{Z}_{\sigma_1}^{(1)}| + \cdots + |\hat{Z}_{\sigma_{[c_1 a_N]}}^{([c_1 a_N])}| \leq [2 c_1 a_N]\}.$$

Since $|\bar{Z}_{\sigma_1}| \leq N_{\sigma_1}^{\bar{Z}}$ $P_0^{N-d\alpha}$-a.s. [cf. (2.13), (2.16)], we have $E_0^{N-d\alpha}[|\bar{Z}_{\sigma_1}|] \leq E_0^{N-d\alpha}[N_{\sigma_1}^{\bar{Z}}] = \frac{1}{d}$ [cf. (2.20)]. Hence, by the law of large numbers,

$$(3.37) \qquad P_0^{N-d\alpha}[A(c_1 a_N)] \overset{N \to \infty}{\longrightarrow} 1.$$

For the probability in (3.35), we obtain the following lower bound using independence of $\{\bar{Z}^{(i)}, \hat{Z}^{(i)}, \sigma_i\}_{i \geq 1}$ and $N^{-d\alpha} > 0$, for $N \geq c(c_1, \varepsilon)$:

$$P_0^{N^{-d\alpha}}\left[\sum_{i=1}^{[c_1 a_N]}(\sigma_i + H_{-\hat{Z}_{\sigma_i}^{(i)}}^{\bar{Z}^{(i)}}) \leq a_N^2(\log N)^\varepsilon\right]$$

$$\overset{(3.36)}{\geq} P_0^{N^{-d\alpha}}\left[\sum_{i=1}^{[c_1 a_N]} H_{-\hat{Z}_{\sigma_i}^{(i)}}^{\bar{Z}^{(i)}} \leq \tfrac{1}{2}a_N^2(\log N)^\varepsilon, A([c_1 a_N])\right]$$

$$\overset{(\text{indep.}, \; N^{-d\alpha}>0)}{\geq} \sum_{j_1=-[2c_1 a_N]}^{[2c_1 a_N]} \cdots \sum_{j_{[c_1 a_N]}=-[2c_1 a_N]}^{[2c_1 a_N]}$$

$$\times P_0^{N^{-d\alpha}}\left[\sum_{i=1}^{[c_1 a_N]} H_{-|j_i|}^{\bar{Z}^{(i)}} \leq \tfrac{1}{2}a_N^2(\log N)^\varepsilon\right]$$

$$\times P_0^{N^{-d\alpha}}[A(c_1 a_N), \hat{Z}_{\sigma_i}^{(i)} = j_i,$$

$$i = 1, \ldots, [c_1 a_N]].$$

(3.38)

By the simple Markov property and the fact that the increments of $\bar{Z}$ are independent and identically distributed, we have the following equality in distribution:

$$\sum_{i=1}^{[c_1 a_N]} H_{-|j_i|}^{\bar{Z}^{(i)}} \overset{(\text{dist.})}{=} H_{-|j_1|-\cdots-|j_{[c_1 a_N]}|}^{\bar{Z}}. \tag{3.39}$$

For the $j_i$'s summed over in (3.38) [recall the definition of $A(c_1 a_N)$ in (3.36)], we have $-|j_1| - \cdots - |j_{[c_1 a_N]}| \geq -[2c_1 a_N]$, so (3.39) implies that, for such $j_i$'s,

$$P_0^{N^{-d\alpha}}\left[\sum_{i=1}^{[c_1 a_N]} H_{-|j_i|}^{\bar{Z}^{(i)}} \leq \tfrac{1}{2}a_N^2(\log N)^\varepsilon\right] \geq P_0^{N^{-d\alpha}}[H_{-[2c_1 a_N]}^{\bar{Z}} \leq \tfrac{1}{2}a_N^2(\log N)^\varepsilon].$$

Substituted into (3.38), this yields

$$P_0^{N^{-d\alpha}}\left[\sum_{i=1}^{[c_1 a_N]}(\sigma_i + H_{-\hat{Z}_{\sigma_i}^{(i)}}^{\bar{Z}^{(i)}}) \leq a_N^2(\log N)^\varepsilon\right]$$

$$\geq P_0^{N^{-d\alpha}}[H_{-[2c_1 a_N]}^{\bar{Z}} \leq \tfrac{1}{2}a_N^2(\log N)^\varepsilon] P_0^{N^{-d\alpha}}[A(c_1 a_N)].$$

(3.40)

Since we already know (3.37), the proof of (3.35), and hence of (3.3), will be complete once it is shown that

$$P_0^{N^{-d\alpha}}[H_{-[2c_1 a_N]}^{\bar{Z}} \leq \tfrac{1}{2}a_N^2(\log N)^\varepsilon] \overset{N\to\infty}{\longrightarrow} 1. \tag{3.41}$$

For $\pi_{\mathbb{Z}}(X)$, the $\mathbb{Z}$-projection of the discrete-time random walk $X$, any $v \in \mathbb{Z}$ and $s, t \geq 0$, we have [cf. (2.13)]

$$\{H_{-v}^{\pi_{\mathbb{Z}}(X)} \leq s, N_t^{\bar{X}} \geq s\} \subseteq \{H_{-v}^{\pi_{\mathbb{Z}}(\bar{X})} \leq t\} = \{H_{-v}^{\bar{Z}} \leq t\}.$$

By this last observation, applied with $t_N = \frac{1}{2} a_N^2 (\log N)^\varepsilon$, $s_N = \frac{d+1}{ed} t_N$, $v = [2c_1 a_N]$, we see with (3.5) and (2.20) that instead of (3.41) it suffices to show that

$$(3.42) \qquad P_0^{N^{-d\alpha}} \left[ H_{-[2c_1 a_N]}^{\pi_{\mathbb{Z}}(X)} \leq \frac{d+1}{2ed} a_N^2 (\log N)^\varepsilon \right] \overset{N \to \infty}{\longrightarrow} 1.$$

With (2.27) of Lemma 2.1, applied with $T = H_{-[2c_1 a_N]}^{\pi_{\mathbb{Z}}(X)}$, $A = \{H_{-[2c_1 a_N]}^{\pi_{\mathbb{Z}}(X)} \leq \frac{d+1}{2ed} a_N^2 (\log N)^\varepsilon\}$ and $b = -[2c_1 a_N]$, we can bound the probability in (3.42) from below by

$$(1 - N^{-d\alpha})^{ca_N} (1 - N^{-2d\alpha})^{ca_N^2 (\log N)^\varepsilon} P_0^0 [H_{-[ca_N]}^{\pi_{\mathbb{Z}}(X)} \leq c' a_N^2 (\log N)^\varepsilon].$$

Since $\alpha > 1$, the factor before the above probability tends to 1 as $N \to \infty$ [cf. (3.32)], while the last probability tends to 1 by the invariance principle. This shows (3.42), hence (3.41), and thus completes the proof of (3.3).

*Case $\alpha \leq 1$.* We claim that in order to prove (3.4), it suffices to show that for some constant $c_2(c_1) > 0$ and $N \geq c(c_1)$, with $a_N$ defined in (3.32),

$$(3.43) \qquad P_0^{N^{-d\alpha}} [\bar{X}([0, N^{3d})) \supseteq \mathbb{T}_N^d \times \{0\}] \geq e^{-c_2 N^{-d\alpha} a_N}$$

(recall our convention concerning constants from the end of the Introduction). Indeed, suppose that (3.43) holds true. Then observe that, on the event $\{\bar{T}^{\text{disc}} \geq e^{c_0 N^{-d\alpha} a_N}\}$, $\bar{X}$ does not cover $\mathbb{T}_N^d \times \{\bar{Z}_{nN^{3d}}\}$ during the time interval $[nN^{3d}, (n+1)N^{3d})$ for $0 \leq n \leq [N^{-3d} e^{c_0 N^{-d\alpha} a_N}] - 1$, $n \geq 1$, for covering of a slice of $E$ results in the disconnection of $E$. We thus apply the simple Markov property inductively at times

$$\{nN^{3d} : n = [N^{-3d} e^{c_0 N^{-d\alpha} a_N}] - 1, \ldots, 2, 1\},$$

and obtain

$$P_0^{N^{-d\alpha}} [\bar{T}_N^{\text{disc}} \geq e^{c_0 N^{-d\alpha} a_N}]$$

$$\leq P_0^{N^{-d\alpha}} \left[ \bigcap_{n=0}^{[N^{-3d} e^{c_0 N^{-d\alpha} a_N}] - 1} \bar{\theta}_{nN^{3d}}^{-1} \{\bar{X}([0, N^{3d})) \not\supseteq \mathbb{T}_N^d \times \{\bar{Z}_{nN^{3d}}\}\} \right]$$

$$\overset{\text{(Markov, transl. inv.)}}{=} P_0^{N^{-d\alpha}} [\bar{X}([0, N^{3d})) \not\supseteq \mathbb{T}_N^d \times \{0\}]^{[N^{-3d} e^{(c_0 N^{-d\alpha} a_N})]}$$

$$\overset{(3.43)}{\leq} \exp\{-ce^{(c_0 - c_2) N^{-d\alpha} a_N} N^{-3d}\}$$

$$\overset{(3.32)}{=} \exp\{-ce^{(c_0 - c_2) N^{d(1-\alpha)} (\log N)^2} N^{-3d}\},$$

so the proof of (3.4) is complete by the fact that $\alpha \leq 1$, provided we choose $c_0 > c_2(c_1)$. It thus remains to establish the estimate (3.43). To this end, we observe that

$$
\begin{aligned}
(3.44) \quad & P_0^{N^{-d\alpha}}[\bar{X}([0, N^{3d})) \supseteq \mathbb{T}_N^d \times \{0\}] \\
& \geq P_0^{N^{-d\alpha}}[\bar{X}([0, \infty)) \supseteq \mathbb{T}_N^d \times \{0\}] \\
& \quad - P_0^{N^{-d\alpha}}[\bar{X}([N^{3d}, \infty)) \cap \mathbb{T}_N^d \times \{0\} \neq \varnothing].
\end{aligned}
$$

Standard large deviations estimates allow us to bound the second probability on the right-hand side. Observe that independence of $N^{\bar{Z}}$ and $Z$ [cf. (2.14)] implies with Fubini's theorem that

$$
\begin{aligned}
& P_0^{N^{-d\alpha}}[\bar{X}([N^{3d}, \infty)) \cap \mathbb{T}_N^d \times \{0\} \neq \varnothing] \\
& = P_0^{N^{-d\alpha}}[\text{for some } t \geq N^{3d}, \bar{Z}_t = 0] \\
& \overset{(2.13)}{=} P_0^{N^{-d\alpha}}[\text{for some } k \geq N_{N^{3d}}^{\bar{Z}}, Z_k = 0] \\
& \overset{\text{(Fubini)}}{=} E_0^{N^{-d\alpha}}[P_0^{N^{-d\alpha}}[\text{for some } k \geq n, Z_k = 0]|_{n = N_{N^{3d}}^{\bar{Z}}}] \\
& \leq E_0^{N^{-d\alpha}}\left[ \sum_{k \geq N_{N^{3d}}^{\bar{Z}}} P_0^{N^{-d\alpha}}[Z_k - N^{-d\alpha}k < -\tfrac{1}{2}N^{-d\alpha}k] \right].
\end{aligned}
$$

Now observe that $(Z_n - \Delta n)_{n \geq 0}$ is a $P_0^{\Delta}$-martingale with increments bounded by $1 + \Delta \leq 2$ [cf. (2.16)]. By Azuma's inequality (see, e.g., [2], page 85), the expression in the last sum is therefore bounded from above by $\exp\{-cN^{-2d\alpha}k\}$. This yields

$$
\begin{aligned}
& P_0^{N^{-d\alpha}}[\bar{X}([N^{3d}, \infty)) \cap \mathbb{T}_N^d \times \{0\} \neq \varnothing] \\
& \leq E_0^{N^{-d\alpha}}\left[ \sum_{k \geq N_{N^{3d}}^{\bar{Z}}} e^{-cN^{-2d\alpha}k} \right] \\
& = \frac{1}{1 - e^{-cN^{-2d\alpha}}} E_0^{N^{-d\alpha}}[e^{-cN^{-2d\alpha} N_{N^{3d}}^{\bar{Z}}}] \\
& \overset{(N \geq c'(\alpha))}{\leq} cN^{2d\alpha} E_0^{N^{-d\alpha}}[e^{-cN^{-2d\alpha} N_{N^{3d}}^{\bar{Z}}}] \\
& \overset{(2.20)}{=} cN^{2d\alpha} \exp\left\{ -\frac{1}{d} N^{3d}(1 - e^{-cN^{-2d\alpha}}) \right\} \\
& \overset{(N \geq c'(\alpha))}{\leq} \exp\{-cN^{3d}N^{-2d\alpha}\} \overset{(\alpha \leq 1)}{\leq} \exp\{-cN^d\}.
\end{aligned}
$$

Inserting this last estimate into (3.44), we see that in fact (3.43) will follow from

$$(3.45) \qquad P_0^{N-d\alpha}[\bar{X}([0,\infty)) \supseteq \mathbb{T}_N^d \times \{0\}] \geq c e^{-(1/2)c_2 N^{-d\alpha} a_N}.$$

By (3.26), we have

$$
\begin{aligned}
(3.46) \qquad & P_0^{N-d\alpha}[\bar{X}([0,\infty)) \supseteq \mathbb{T}_N^d \times \{0\}] \\
& \geq P_0^{N-d\alpha}[\bar{S}_{[c_1 a_N]} < \infty, \bar{X}([0, \bar{S}_{[c_1 a_N]})) \supseteq \mathbb{T}_N^d \times \{0\}] \\
& = P_0^{N-d\alpha}[\bar{S}_{[c_1 a_N]} < \infty] P_0^{N-d\alpha}[C_{\mathbb{T}_N^d}^V \leq [c_1 a_N] | \bar{S}_{[c_1 a_N]} < \infty] \\
& \overset{((3.26),(3.32))}{\geq} c P_0^{N-d\alpha}[\bar{S}_{[c_1 a_N]} < \infty].
\end{aligned}
$$

With the help of (3.12), we obtain, with the same arguments as in (3.38), (3.39) and (3.40) with $A(c_1 a_N)$ replaced by $\{\sum_{i=1}^{[c_1 a_N]} |\hat{Z}_{\sigma_i}^{(i)}| \leq [2c_1 a_N]\}$,

$$
\begin{aligned}
(3.47) \qquad & P_0^{N-d\alpha}[\bar{S}_{[c_1 a_N]} < \infty] \\
& \overset{(3.12)}{=} P_0^{N-d\alpha}\left[ \sum_{i=1}^{[c_1 a_N]} H_{-\hat{Z}_{\sigma_i}^{(i)}}^{\bar{Z}^{(i)}} < \infty \right] \\
& \geq P_0^{N-d\alpha}[H_{-[2c_1 a_N]}^{\bar{Z}} < \infty] P_0^{N-d\alpha}\left[ \sum_{i=1}^{[c_1 a_N]} |\hat{Z}_{\sigma_i}^{(i)}| \leq [2c_1 a_N] \right] \\
& = \left( \frac{1 - N^{-d\alpha}}{1 + N^{-d\alpha}} \right)^{[2c_1 a_N]} P_0^{N-d\alpha}\left[ \sum_{i=1}^{[c_1 a_N]} |\hat{Z}_{\sigma_i}^{(i)}| \leq [2c_1 a_N] \right].
\end{aligned}
$$

The factor in front of the probability on the right-hand side is bounded from below by $e^{-c(c_1)N^{-d\alpha} a_N}$, while the probability tends to 1 as $N \to \infty$, again by the estimate $E_0^{N-d\alpha}[|\bar{Z}_{\sigma_1}|] \leq E_0^{N-d\alpha}[N_{\sigma_1}^{\bar{Z}}] = \frac{1}{d}$ and the law of large numbers. Therefore, (3.46) and (3.47) together show (3.45) for a suitably chosen constant $c_2(c_1) > 0$. Hence, the proof of (3.4) and thus of Theorem 3.2 is complete. $\square$

## 4. Lower bounds: Reduction to large deviations.
The goal of this section is to prove Theorem 1.2 reducing the problem of finding a lower bound on $T_N^{\text{disc}}$ to a large deviations estimate of the form (1.10). As a preliminary step toward this reduction, we prove the following geometric lemma in the spirit of Dembo and Sznitman [4], where we refer to (1.6) for our notion of $\kappa$-disconnection:

LEMMA 4.1 [$d \geq 1$, $\alpha > 0$, $\kappa \in (0, \frac{1}{2})$]. *There is a constant $c(\alpha, \kappa)$ such that for all $N \geq c(\alpha, \kappa)$, whenever $K \subseteq E$ disconnects $E$, there is an $x_* \in E$* such that $K$ $\kappa$-disconnects $x_* + B(\alpha)$; cf. (1.7). (We refer to the end of the *Introduction* for our convention concerning constants.)

PROOF. We follow the argument contained in the proof of Lemma 2.4 in Dembo and Sznitman [4]. Assuming that $K$ disconnects $E$, we refer as *Top* to the connected component of $E \setminus K$ containing $\mathbb{T}_N^d \times [M, \infty)$ for large $M \geq 1$. We can then define the function

$$t : E \longrightarrow \mathbb{R}_+$$
$$x \mapsto \frac{|Top \cap (x + B(\alpha))|}{|B(\alpha)|}.$$

The function $t$ takes the value 0 for $x = (u, v) \in E$ with $v \in \mathbb{Z}$ a large negative number and the value 1 for $v$ a large positive number. Moreover, for $x = (u, v)$, $x' = (u, v') \in E$ such that $|v - v'| = 1$ we have (with $\Delta$ denoting symmetric difference)

$$|t(x) - t(x')| \leq \frac{|(x + B(\alpha)) \Delta (x' + B(\alpha))|}{|B(\alpha)|} \leq \frac{cN^d}{N^{d + d\alpha \wedge 1}} = \frac{c}{N^{d\alpha \wedge 1}}.$$

Using these last two observations on $t$, we see that, for $N \geq c(\alpha, \kappa)$, there is at least one $x_* \in E$ satisfying

$$\left| t(x_*) - \frac{1}{2} \right| \leq \frac{c}{N^{d\alpha \wedge 1}} \leq \frac{1}{2} - \kappa,$$

which can be restated as

(4.1)  $$\kappa |B(\alpha)| \leq |Top \cap (x_* + B(\alpha))| \leq (1 - \kappa)|B(\alpha)|.$$

If we set $I = Top \cap (x_* + B(\alpha))$, then $\partial_{(x_* + B(\alpha))}(I) \subseteq K$ (since $K$ disconnects $E$), so that the proof is complete with (4.1). □

PROOF OF THEOREM 1.2. We claim that it suffices to prove the following two estimates on $P_0^{N^{-d\alpha}}[T_N^{\mathrm{disc}} \leq t]$, valid for any $t \geq 1$, $\xi \in (0, f(\alpha, \beta))$ [for $\alpha, \beta > 0$ and $f$ as in (1.10)] and $N \geq c(\alpha, \beta, \xi)$:

(4.2)  $$P_0^{N^{-d\alpha}}[T_N^{\mathrm{disc}} \leq t] \leq cN^d(t + N)(e^{-N^\xi} + e^{-c'N^{\beta + (d\alpha \wedge 1)}t^{-1/2}})$$

and

(4.3)  $$P_0^{N^{-d\alpha}}[T_N^{\mathrm{disc}} \leq t] \leq cN^d(t + N)(e^{-N^\xi} + e^{-c'N^{\beta - (d\alpha - 1)_+}}).$$

Indeed, suppose that (4.2) and (4.3) both hold. In order to deduce (1.11), we then choose any $\alpha > 1$, $0 < \varepsilon < 2d$ such that $\beta = d - 1 - \frac{\varepsilon}{4} > 0$ (note $d \geq 2$) and $\xi \in (0, f(\alpha, \beta))$ (which is possible by the assumption on $f$). With $t = N^{2d - \varepsilon}$, (4.2) then yields, for $N \geq c(\alpha, \beta, \xi, \varepsilon)$,

$$P_0^{N^{-d\alpha}}[T_N^{\mathrm{disc}} \leq N^{2d - \varepsilon}] \leq cN^{3d - \varepsilon}(e^{-N^\xi} + e^{-c'N^{\varepsilon/4}}),$$

and hence shows (1.11).

On the other hand, choosing $t = \exp\{N^\mu\}$, $\mu > 0$, in (4.3), we have, for any $\alpha, \beta > 0$, $\xi \in (0, f(\alpha, \beta))$ and $N \geq c(\alpha, \beta, \xi, \mu)$,

$$
\begin{aligned}
(4.4) \quad & P_0^{N^{-d\alpha}}[T_N^{\mathrm{disc}} \leq \exp\{N^\mu\}] \\
& \leq cN^d(\exp\{N^\mu - N^\xi\} + \exp\{N^\mu - c'N^{\beta - (d\alpha - 1)_+}\}).
\end{aligned}
$$

The right-hand side of (4.4) tends to 0 as $N \to \infty$ for $\alpha, \beta, \xi$ as above, provided $\beta > (d\alpha - 1)_+$ and $\mu < \xi \wedge (\beta - (d\alpha - 1)_+)$. We thus obtain (1.12) by optimizing over $\beta$ and $\xi$ in (4.4).

It therefore remains to establish (4.2) and (4.3). To this end, we apply the geometric Lemma 4.1, noting that, up to time $t$, only sets $(u, v) + B(\alpha)$ [in the notation of (1.7)] with $|v| \leq t + N^{d\alpha \wedge 1}$ can be entered by the discrete-time random walk, and thus deduce that, for $N \geq c(\alpha)$,

$$
\begin{aligned}
(4.5) \quad & P_0^{N^{-d\alpha}}[T_N^{\mathrm{disc}} \leq t] \\
& \leq cN^d(t + N) \sup_{x \in E} P_0^{N^{-d\alpha}}[X([0, [t]])\ \tfrac{1}{3}\text{-disconnects } x + B(\alpha)].
\end{aligned}
$$

For the first return time $\mathcal{R}_1^x$, defined as $\mathcal{R}_1^x = H_{S_{2[N^{d\alpha \wedge 1}]}}^X$ [cf. (2.22)], one has

$$
\begin{aligned}
& \{X([0, [t]])\ \tfrac{1}{3}\text{-disconnects } x + B(\alpha)\} \\
& \subseteq \theta_{\mathcal{R}_1^x}^{-1}\{X([0, [t]])\ \tfrac{1}{3}\text{-disconnects } x + B(\alpha)\}.
\end{aligned}
$$

Applying the strong Markov property at time $\mathcal{R}_1^x$ and using translation invariance, we thus obtain that [cf. (1.8)]

$$
\begin{aligned}
& P_0^{N^{-d\alpha}}[X([0, [t]])\ \tfrac{1}{3}\text{-disconnects } x + B(\alpha)] \\
& \leq \sup_{x \in S_{2[N^{d\alpha \wedge 1}]}} P_x^{N^{-d\alpha}}[X([0, [t]])\ \tfrac{1}{3}\text{-disconnects } B(\alpha)] \\
& = \sup_{x \in S_{2[N^{d\alpha \wedge 1}]}} P_x^{N^{-d\alpha}}[U_{B(\alpha)} \leq t].
\end{aligned}
$$

Inserted into (4.5), this yields

$$
(4.6) \quad P_0^{N^{-d\alpha}}[T_N^{\mathrm{disc}} \leq t] \leq cN^d(t + N) \sup_{x \in S_{2[N^{d\alpha \wedge 1}]}} P_x^{N^{-d\alpha}}[U_{B(\alpha)} \leq t].
$$

We then observe that, for any $x \in S_{2[N^{d\alpha \wedge 1}]}$,

$$
\begin{aligned}
(4.7) \quad & P_x^{N^{-d\alpha}}[U_{B(\alpha)} \leq t] \\
& \leq P_x^{N^{-d\alpha}}[U_{B(\alpha)} < \mathcal{D}_{[N^\beta]}] + P_x^{N^{-d\alpha}}[\mathcal{R}_{[N^\beta]} \leq U_{B(\alpha)} \leq t] \\
& \overset{\mathrm{(def.)}}{=} P_1 + P_2.
\end{aligned}
$$

By definition of $U_{B(\alpha)}$ we know that, on the event $\{U_{B(\alpha)} < \infty\}$, $\pi_{\mathbb{Z}}(X_{U_{B(\alpha)}} - x) \leq c[N^{d\alpha \wedge 1}]$, $P_x^{N-d\alpha}$-a.s., for $x \in S_{2[N^{d\alpha \wedge 1}]}$. We can thus apply (2.28) of Lemma 2.1 with $A = \{U_{B(\alpha)} < \mathcal{D}_{[N^\beta]}\}$, $T = U_{B(\alpha)}$ and $b' = c[N^{d\alpha \wedge 1}]$ and obtain, for $P_1$ in (4.7),

$$
\begin{aligned}
P_1 &\overset{(2.28)}{\leq} (1 + N^{-d\alpha})^{c[N^{d\alpha \wedge 1}]} P_x^0[U_{B(\alpha)} < \mathcal{D}_{[N^\beta]}] \\
&\leq c P_x^0[U_{B(\alpha)} < \mathcal{D}_{[N^\beta]}] \\
&\overset{(1.10)}{\leq} e^{-N^\xi},
\end{aligned}
$$
(4.8)

for any $\xi \in (0, f(\alpha, \beta))$ and all $N \geq c(\alpha, \beta, \xi)$. Turning to $P_2$ in (4.7), we apply (2.28) of Lemma 2.1 with $A = \{\mathcal{R}_{[N^\beta]} \leq t\}$, $T = \mathcal{R}_{[N^\beta]}$ and $b' = c[N^{d\alpha \wedge 1}]$, and obtain

$$
\begin{aligned}
P_2 &\leq P_x^{N-d\alpha}[\mathcal{R}_{[N^\beta]} \leq t] \leq (1 + N^{-d\alpha})^{c[N^{d\alpha \wedge 1}]} P_x^0[\mathcal{R}_{[N^\beta]} \leq t] \\
&\leq c P_x^0[\mathcal{R}_{[N^\beta]} \leq t].
\end{aligned}
$$

For this last probability, we make the observation that, under $P_x^0$, $\mathcal{R}_{[N^\beta]} - \mathcal{D}_1$ ($\leq \mathcal{R}_{[N^\beta]}$) is distributed as the sum of at least $[cN^{d\alpha \wedge 1} N^\beta]$ independent random variables, all of which are distributed as the hitting time of 1 for the unbiased simple random walk $\pi_{\mathbb{Z}}(X)$ [cf. (2.1)] starting at the origin with geometric delay of constant parameter $\frac{1}{d+1}$. Applying an elementary estimate on one-dimensional simple random walk for the second inequality (cf. Durrett [7], Chapter 3, (3.4)), we deduce that, for $t \geq 1$,

$$
\begin{aligned}
P_2 &\leq c P_0^0[H_1^{\pi_{\mathbb{Z}}(X)} \leq t]^{cN^{\beta + (d\alpha \wedge 1)}} \leq c(1 - c't^{-1/2})^{c'N^{\beta + (d\alpha \wedge 1)}} \\
&\leq c \exp\{-c'N^{\beta + (d\alpha \wedge 1)} t^{-1/2}\}.
\end{aligned}
$$

Together with (4.8), (4.7) and (4.6), this yields (4.2).

In order to obtain (4.3), we use the following different method for estimating $P_2$ in (4.7): We let $A^-$ be the event that the random walk $X$ first exits $S_{4[N^{d\alpha \wedge 1}]}$ into the negative direction, that is,

$$
A^- = \{\pi_{\mathbb{Z}}(X_{\mathcal{D}_1}) < 0\} \in \mathcal{F}_{\mathcal{D}_1}.
$$

One then has

$$
\begin{aligned}
P_2 &\leq \sup_{x \in S_{2[N^{d\alpha \wedge 1}]}} P_x^{N-d\alpha}[\mathcal{R}_{[N^\beta]} < \infty] \\
&= \sup_{x \in S_{2[N^{d\alpha \wedge 1}]}} (P_x^{N-d\alpha}[\mathcal{R}_{[N^\beta]} < \infty, A^-] + P_x^{N-d\alpha}[\mathcal{R}_{[N^\beta]} < \infty, (A^-)^c]).
\end{aligned}
$$
(4.9)

We now apply the strong Markov property at the times $\mathcal{D}_1$ and $\mathcal{R}_2$ and use translation invariance to infer from (4.9) that, for $N^\beta \geq 2$,

$$
\begin{aligned}
(4.10) \quad P_2 &\leq \sup_{x \in S_{2[N^{d\alpha \wedge 1}]}} P_x^{N-d\alpha}[\mathcal{R}_{[N^\beta]-1} < \infty] \\
&\times \sup_{x \in S_{2[N^{d\alpha \wedge 1}]}} (P_x^{N-d\alpha}[A^-] \\
&\qquad + P_x^{N-d\alpha}[(A^-)^c] P_0^{N-d\alpha}[H_{-c[N^{d\alpha \wedge 1}]}^{\pi_{\mathbb{Z}}(X)} < \infty]).
\end{aligned}
$$

Next, we apply the estimate (2.28) of Lemma 2.1 with $T = \mathcal{D}_1$, $A = A^-$ and $b' = -2[N^{d\alpha \wedge 1}]$, then the invariance principle for one-dimensional simple random walk, and obtain, for any $x \in S_{2[N^{d\alpha \wedge 1}]}$,

$$
(4.11) \qquad P_x^{N-d\alpha}[A^-] \overset{(2.28)}{\leq} P_x^0[A^-] \overset{\text{(inv. princ.)}}{\leq} (1 - c_3), \qquad c_3 > 0.
$$

Moreover, since the projection $\pi_{\mathbb{Z}}(X)$ of $X$ on $\mathbb{Z}$ is a one-dimensional random walk with drift $\frac{N^{-d\alpha}}{d+1}$ and geometric delay of constant parameter $\frac{1}{d+1}$, standard estimates on one-dimensional biased random walk imply

$$
\begin{aligned}
(4.12) \quad P_0^{N-d\alpha}[H_{-c[N^{d\alpha \wedge 1}]}^{\pi_{\mathbb{Z}}(X)} < \infty] &\leq \left( \frac{1 - N^{-d\alpha}(d+1)^{-1}}{1 + N^{-d\alpha}(d+1)^{-1}} \right)^{c[N^{d\alpha \wedge 1}]} \\
&\leq e^{-cN^{-d\alpha}[N^{d\alpha \wedge 1}]}.
\end{aligned}
$$

Inserting (4.11) and (4.12) into (4.10) and using induction, we deduce

$$
(4.13) \qquad P_2 \leq (1 - c_3 + c_3 e^{-cN^{-d\alpha}[N^{d\alpha \wedge 1}]})^{[N^\beta]-1}.
$$

Note that $N^{-d\alpha}[N^{d\alpha \wedge 1}] \leq 1$. If $d\alpha > 1$, then the right-hand side of (4.13) is bounded from above by $(1 - cN^{-d\alpha}[N^{d\alpha \wedge 1}])^{[N^\beta]-1} \leq e^{-cN^{-d\alpha}[N^{d\alpha \wedge 1}]N^\beta}$, while if $d\alpha \leq 1$, the right-hand side of (4.13) is bounded by $e^{-cN^\beta}$. In any case, we infer from (4.13) that

$$
P_2 \leq e^{-cN^{-d\alpha}N^{\beta+(d\alpha \wedge 1)}} = e^{-cN^{\beta-(d\alpha-1)_+}}.
$$

Together with (4.8), (4.7) and (4.6) this yields (4.3) and completes the proof of Theorem 1.2. $\quad \square$

## 5. More geometric lemmas.

The purpose of this section is to prove several geometric lemmas needed for the derivation of the large deviations estimate (1.10) in Theorem 1.2. The general purpose of these geometric results is to impose restrictions on a set $K$ $\frac{1}{3}$-disconnecting $B(\alpha)$. This will enable us to obtain an upper bound on the probability appearing in (1.10), when choosing $K = X([0, \mathcal{D}_{[N^\beta]}])$.

Throughout this and the next section, we consider the scales $L$ and $l$, defined as

$$(5.1) \qquad l = [N^\gamma], \qquad L = [N^{\gamma'}] \qquad \text{for } 0 < \gamma < \gamma' \wedge d\alpha, 0 < \gamma' < 1.$$

The crucial geometric estimates come in Lemma 5.3 and its modification Lemma 5.4. These geometric results, in the spirit of Dembo and Sznitman [4], require as key ingredient an isoperimetric inequality of Deuschel and Pisztora [6]; see Lemma 5.2. In rough terms, Lemmas 5.3 and 5.5 show that for any set $K$ disconnecting $C(L)$ or $B(\alpha)$ for $d\alpha < 1$ [cf. (1.7), (2.7)], one can find a whole "surface" of subcubes of $C(L)$ or $B(\alpha)$ such that the set $K$ occupies a "surface" of points inside every one of these subcubes. More precisely, it is shown that there exist subcubes $(C_x(l))_{x \in \mathcal{E}}$ [cf. (2.8)] of $C(L)$, respectively of $B(\alpha)$, with the following properties: for one of the projections $\pi_*$ on the $d$-dimensional hyperplanes, the projected set of base-points $\pi_*(\mathcal{E})$ is arranged on a subgrid of side-length $l$ and is substantially large. In the case of $C(L)$, this set of points occupies at least a constant fraction of the volume of the projected subgrid of $C(L)$. Moreover, for one of the projections $\pi_{**}$ (possibly different from $\pi_*$), the $\pi_{**}$-projection of the disconnecting set $K$ intersected with any subcube $C_x(l)$, $x \in \mathcal{E}$, contains at least $cl^d$ points, that is, at least a constant fraction of the volume of $\pi_{**}(C_x(l))$ (see Figure 3 for an illustration of the idea).

The first lemma in this section allows to propagate disconnection of the $|\cdot|_\infty$-ball $B_\infty(0, [N/4])$ to a smaller scale of size $L$, in the sense that, for any set $K$ $\frac{1}{3}$-disconnecting $B_\infty(0, [N/4])$, one can find a sub-box $C_{x_*}(L)$ of $B_\infty(0, [N/4])$ which is $\frac{1}{4}$-disconnected by $K$ [cf. (1.6)]. This result will prove useful for the case $B(\alpha) = B_\infty(0, [N/4])$ (i.e., if $d\alpha \geq 1$), where we use an upper bound on the number of excursions between $C_{x_*}(L)$ and $(C_{x_*}(L)^{(L)})^c$ performed by the random walk $X$ until time $\mathcal{D}_{[N^\beta]}$. We refer to the end of the Introduction for our convention concerning constants.

LEMMA 5.1 ($d \geq 1$, $\gamma' \in (0,1)$, $L = [N^{\gamma'}]$, $N \geq 1$). *There is a constant $c(\gamma') > 0$ such that for all $N \geq c(\gamma')$, whenever $K \subseteq B_\infty(0, [N/4])$ $\frac{1}{3}$-disconnects $B_\infty(0, [N/4])$, there is an $x_* \in B_\infty(0, [N/4])$ such that $K$ $\frac{1}{4}$-disconnects $C_{x_*}(L) \subseteq B_\infty(0, [N/4])$.*

PROOF. Since $K$ $\frac{1}{3}$-disconnects $B_\infty(0, [N/4])$ [cf. (1.6)], there is a set $I \subseteq B_\infty(0, [N/4])$ satisfying $\frac{1}{3}|B_\infty(0, [N/4])| \leq |I| \leq \frac{2}{3}|B_\infty(0, [N/4])|$ and $\partial_{B_\infty(0, [N/4])}(I) \subseteq K$. We want to find a point $x_* \in E$ such that $C_{x_*}(L) \subseteq B_\infty(0, [N/4])$ and

$$(5.2) \qquad \frac{1}{4}|C(L)| \leq |C_{x_*}(L) \cap I| \leq \frac{3}{4}|C(L)|.$$

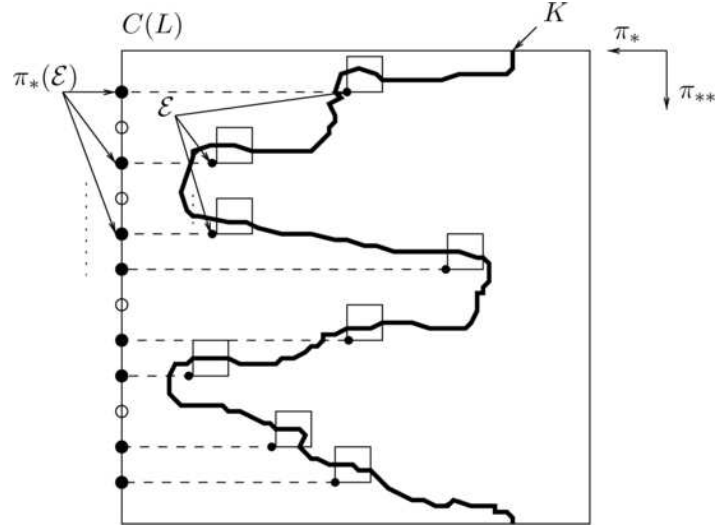

Fig. 3. *An illustration of the crucial geometric Lemma 5.3. The figure shows the set $C(L)$, disconnected by $K \subseteq C(L)$. The small boxes are the collection of subcubes $(C_x(l))_{x \in \mathcal{E}}$. The circles on the left are the points on the projected subgrid of side-length $l$, a large number of which (the filled ones) are occupied by the projected set $\pi_*(\mathcal{E})$ of base-points $\mathcal{E}$ [cf. (5.12), (5.13)]. In every subcube, the set $K$ occupies a surface of a significant number of points, in the sense of (5.14).*

To this end, we introduce the subgrid $\mathcal{B}_L \subseteq B_\infty(0, [N/4])^{(-L)}$ of side-length $L$, defined as

$$(5.3) \qquad \mathcal{B}_L = B_\infty(0, [N/4])^{(-L)} \cap \pi_E([-[N/4], [N/4]]^{d+1} \cap L\mathbb{Z}^{d+1})$$

[cf. (2.3)].

The boxes $(C_x(L))_{x \in \mathcal{B}_L}$ [see (2.7), (2.8)] are disjoint subsets of $B_\infty(0, [N/4])$, and their union covers all but at most $cN^d L$ points of $B_\infty(0, [N/4])$. Hence, we have

$$(5.4) \qquad \sum_{x \in \mathcal{B}_L} |I \cap C_x(L)| \leq |I| \leq \sum_{x \in \mathcal{B}_L} |I \cap C_x(L)| + cN^d L,$$

$$(5.5) \qquad |B_\infty(0, [N/4])| - cN^d L \leq |\mathcal{B}_L||C(L)| \leq |B_\infty(0, [N/4])|.$$

We now claim that, for $N \geq c(\gamma')$, there is at least one $x_1 \in \mathcal{B}_L$ such that

$$(5.6) \qquad |I \cap C_{x_1}(L)| \leq \tfrac{3}{4}|C(L)|.$$

Indeed, otherwise it would follow from the definition of $I$ and the left-hand side inequalities of (5.4) and (5.5) that

$$\tfrac{2}{3}|B_\infty(0, [N/4])| \geq |I| \overset{(5.4)}{>} \tfrac{3}{4}|C(L)||\mathcal{B}_L| \overset{(5.5)}{\geq} \tfrac{3}{4}|B_\infty(0, [N/4])| - cN^d L,$$

which due to the definition of $L$ is impossible for $N \geq c(\gamma')$. Similarly, for $N \geq c(\gamma')$, we can find an $x_2 \in \mathcal{B}_L$ such that

$$(5.7) \qquad \tfrac{1}{4}|C(L)| \leq |I \cap C_{x_2}(L)|,$$

for otherwise the right-hand side inequalities of (5.4) and (5.5) would yield that $\tfrac{1}{3}|B_\infty(0, [N/4])| \leq \tfrac{1}{4}|B_\infty(0, [N/4])| + cN^d L$, thus again leading to a contradiction.

Next, we note that, for any neighbors $x$ and $x'$ in $B_\infty(0, [N/4])$, one has, with $\Delta$ denoting the symmetric difference,

$$(5.8) \qquad \left| \frac{|C_x(L) \cap I|}{|C(L)|} - \frac{|C_{x'}(L) \cap I|}{|C(L)|} \right| \leq \frac{|C_x(L) \Delta C_{x'}(L)|}{|C(L)|} \leq \frac{c}{N^{\gamma'}}.$$

Since both $x_1$ and $x_2$ are in $\mathcal{B}_L \subseteq B_\infty(0, [N/4])^{(-L)}$, we can now choose a nearest-neighbor path $\mathcal{P} = (x_1 = y_1, y_2, \ldots, y_n = x_2)$ from $x_1$ to $x_2$ such that $C_{y_i}(L) \subseteq B_\infty(0, [N/4])$ for all $y_i \in \mathcal{P}$. Consider now the first point $x_* = y_{i_*}$ on $\mathcal{P}$ such that $\tfrac{1}{4}|C(L)| \leq |C_{x_*}(L) \cap I|$, which is well defined thanks to (5.7). If $x_* = y_1$, then by (5.6), $x_*$ satisfies (5.2). If $x_* \neq y_1$, then by (5.8) and choice of $x_*$, one also has

$$\tfrac{1}{4}|C(L)| \leq |C_{x_*}(L) \cap I| \overset{(5.8)}{\leq} |C_{y_{i_*-1}}(L) \cap I| + \frac{c}{N^{\gamma'}}|C(L)|$$

$$< \left( \frac{1}{4} + \frac{c}{N^{\gamma'}} \right)|C(L)|,$$

hence again (5.2) for $N \geq c(\gamma')$. For $N \geq c(\gamma')$, we have thus found an $x_* \in B_\infty(0, [N/4])$ satisfying $\tfrac{1}{4}|C(L)| \leq |C_{x_*}(L) \cap I| \leq \tfrac{3}{4}|C(L)|$ and $C_{x_*}(L) \subseteq B_\infty(0, [N/4])$. Moreover, $\partial_{C_{x_*}(L)}(C_{x_*}(L) \cap I) \subseteq \partial_{B_\infty(0, [N/4])}(I) \subseteq K$. In other words, $K$ $\tfrac{1}{4}$-disconnects $C_{x_*}(L) \subseteq B_\infty(0, [N/4])$. $\square$

The following lemma contains the essential ingredients for the proof of the two main geometric lemmas thereafter.

LEMMA 5.2 $[d \geq 1, \kappa \in (0,1), M \in \{0, \ldots, N-1\}, N \geq 1]$. *Suppose* $A \subseteq [0, M]^{d+1} \subseteq E$. *Then there is an* $i_0 \in \{1, \ldots, d+1\}$ *such that*

$$(5.9) \qquad |A| \leq |\pi_{i_0}(A)|^{(d+1)/d}.$$

*If $A$ in addition satisfies*

$$(5.10) \qquad |A| \leq (1-\kappa)(M+1)^{d+1},$$

*then there is an $i_1 \in \{1, \ldots, d+1\}$ and a constant $c(\kappa) > 0$ such that [cf. (2.6)]*

$$(5.11) \qquad |\pi_{i_1}(\partial_{[0,M]^{d+1}, i_1}(A))| \geq c(\kappa)|A|^{d/(d+1)}.$$

PROOF.   The estimate (5.9) follows, for instance, from a theorem of Loomis and Whitney [10]. The proof of (5.11) can be found in (A.3)–(A.6) in Deuschel and Pisztora [6], page 480.   □

We now come to the main geometric lemma, which provides a necessary criterion for disconnection of the box $C(L)$ [cf. (2.7)]. A schematic illustration of its content can be found in Figure 3.

LEMMA 5.3 ($d \geq 1$, $0 < \gamma < \gamma' < 1$, $l = [N^\gamma]$, $L = [N^{\gamma'}]$, $N \geq 1$).   *For all $N \geq c(\gamma, \gamma')$, whenever $K \subseteq C(L)$ $\frac{1}{4}$-disconnects $C(L)$ [cf. (1.6)], then there exists a set $\mathcal{E} \subseteq C(L)^{(-l)}$ [cf. (2.3)] and projections $\pi_*$ and $\pi_{**} \in \{\pi_1, \ldots, \pi_{d+1}\}$ such that*

$$(5.12) \qquad \pi_*(\mathcal{E}) \subseteq \pi_*(C(L) \cap \pi_E([0, L]^{d+1} \cap l\mathbb{Z}^{d+1})),$$

$$(5.13) \qquad |\pi_*(\mathcal{E})| \geq c' \left( \frac{L}{l} \right)^d,$$

$$(5.14) \qquad \text{for all } x \in \mathcal{E} : |\pi_{**}(K \cap C_x(l))| \geq c'' l^d \qquad [cf. (2.8)].$$

PROOF.   Since $K$ $\frac{1}{4}$-disconnects $C(L)$, there exists a set $I \subseteq C(L)$ satisfying $\frac{1}{4}L^{d+1} \leq |I| \leq \frac{3}{4}L^{d+1}$ and $\partial_{C(L)}(I) \subseteq K$. We introduce here the subgrid $\mathcal{C}_l \subseteq C(L)^{(-l)}$ of side-length $l$, that is,

$$(5.15) \qquad \mathcal{C}_l = C(L)^{(-l)} \cap \pi_E([0, L]^{d+1} \cap l\mathbb{Z}^{d+1}),$$

with sub-boxes $C_x(l)$, $x \in \mathcal{C}_l$. The set $\mathcal{A}$ is then defined as the set of all $x \in \mathcal{C}_l$ whose corresponding box $C_x(l)$ is filled up to more than $\frac{1}{8}$th by $I$:

$$(5.16) \qquad \mathcal{A} = \{x \in \mathcal{C}_l : |C_x(l) \cap I| > \tfrac{1}{8} l^{d+1}\}.$$

Since the disjoint union of the boxes $(C_x(l))_{x \in \mathcal{C}_l}$ contains all but at most $cL^d l$ points of $C(L)$, we have

$$(5.17) \qquad \tfrac{1}{4}L^{d+1} \leq |I| \leq \tfrac{1}{8} l^{d+1} |\mathcal{C}_l \setminus \mathcal{A}| + l^{d+1}|\mathcal{A}| + cL^d l.$$

Using the estimate $|\mathcal{C}_l \setminus \mathcal{A}| \leq |\mathcal{C}_l| \leq (\frac{L}{l})^{d+1}$ and rearranging, we deduce from (5.17) that

$$\left( \frac{1}{8} - c\frac{l}{L} \right) \left( \frac{L}{l} \right)^{d+1} \leq |\mathcal{A}|,$$

so that for $N \geq c(\gamma, \gamma')$,

$$(5.18) \qquad \frac{1}{9}|\mathcal{C}_l| \leq \frac{1}{9} \left( \frac{L}{l} \right)^{d+1} \leq |\mathcal{A}|.$$

In order to apply the isoperimetric inequality (5.11) of Lemma 5.2 with $\mathcal{A}$ and $\mathcal{C}_l$ playing the roles of $A$ and $[0, M]^{d+1}$ for $N \geq c(\gamma, \gamma')$, we need to keep $|\mathcal{A}|$ away from $|\mathcal{C}_l|$. We therefore distinguish two cases, as to whether or not

$$(5.19) \qquad |\mathcal{A}| \leq c_4 |\mathcal{C}_l| \qquad \text{with } c_4^3 = \tfrac{1}{2}(1 + \tfrac{4}{5}).$$

Suppose first that (5.19) holds. Then for $N \geq c(\gamma, \gamma')$, the isoperimetric inequality (5.11), applied on the subgrid $\mathcal{C}_l$, yields an $i \in \{1, \dots, d+1\}$ such that

$$(5.20) \qquad |\pi_i(\partial_{\mathcal{C}_l, i}(\mathcal{A}))| \geq c|\mathcal{A}|^{d/(d+1)} \overset{(5.18)}{\geq} c'\left(\frac{L}{l}\right)^d,$$

where $\partial_{\mathcal{C}_l, i}(\mathcal{A})$ denotes the boundary on the subgrid $\mathcal{C}_l$, defined in analogy with (2.6). In order to construct the set $\mathcal{E}$, we apply the following procedure. Given $w \in \pi_i(\partial_{\mathcal{C}_l, i}(\mathcal{A}))$, we choose an $x' \in \partial_{\mathcal{C}_l, i}(\mathcal{A})$ with $\pi_i(x') = w$. In view of (2.6), at least one of $x' + le_i$ and $x' - le_i$ belongs to $\mathcal{A}$. Without loss of generality, we assume that $x' + le_i \in \mathcal{A}$. We then have $|C_{x'}(l) \cap I| \leq \frac{1}{8}l^{d+1}$ [because $x' \in \mathcal{C}_l \setminus \mathcal{A}$; cf. (5.16)] and $|C_{x'+le_i}(l) \cap I| > \frac{1}{8}l^{d+1}$ (because $x' + le_i \in \mathcal{A}$). Observe that neighboring $x_1, x_2 \in E$ satisfy

$$(5.21) \qquad \left| \frac{|C_{x_1}(l) \cap I|}{l^{d+1}} - \frac{|C_{x_2}(l) \cap I|}{l^{d+1}} \right| \leq \frac{c}{N^\gamma}.$$

Now consider the first point $x = x' + l_* e_i$ on the segment $[x', x' + le_i] = (x', x' + e_i, \dots, x' + le_i)$ satisfying $\frac{1}{8}l^{d+1} < |C_x(l) \cap I|$. By the above observations, this point $x$ is well defined and not equal to $x'$. By (5.21), $x$ then also satisfies

$$(5.22) \qquad \begin{aligned} \frac{l^{d+1}}{8} &< |C_x(l) \cap I| \overset{(5.21)}{\leq} |C_{x'+(l_*-1)e_i}(l) \cap I| + \frac{cl^{d+1}}{N^\gamma} \\ &\leq \left(\frac{1}{8} + \frac{c}{N^\gamma}\right)l^{d+1} \leq \frac{l^{d+1}}{7}, \end{aligned}$$

for $N \geq c(\gamma)$. In addition, one has $\pi_i(x) = \pi_i(x') = w$. This construction thus yields, for any $w \in \pi_i(\partial_{\mathcal{C}_l, i}(\mathcal{A}))$, a point $x \in C(L)^{(-l)}$ [note that $x'$, $x' + le_i \in C(L)^{(-l)}$ and $C(L)^{(-l)}$ is convex], satisfying (5.22) and $\pi_i(x) = w$. We define the set $\mathcal{E}'$ as the set of all such points $x$. Then by construction, we have $\pi_i(\mathcal{E}') = \pi_i(\partial_{\mathcal{C}_l, i}(\mathcal{A}))$; in particular (5.12) holds with $\mathcal{E}'$ in place of $\mathcal{E}$ and $\pi_* = \pi_i$, as does (5.13), by (5.20). For any $x \in \mathcal{E}'$, we apply the isoperimetric inequality (5.11) of Lemma 5.2 with $C_x(l)$ in place of $[0, M]^{d+1}$, $C_x(l) \cap I$ in place of $A$ and $1 - \kappa = \frac{1}{7}$; cf. (5.22). We thus find a $j(x) \in \{1, \dots, d+1\}$ with

$$(5.23) \quad |\pi_{j(x)}(\partial_{C_x(l), j(x)}(C_x(l) \cap I))| \geq c|C_x(l) \cap I|^{d/(d+1)} \overset{(5.22)}{\geq} c'l^d.$$

It follows from the choice of $I$ that $\partial_{C_x(l),j(x)}(C_x(l) \cap I) \subseteq K \cap C_x(l)$, and hence

$$(5.24) \qquad |\pi_{j(x)}(K \cap C_x(l))| \geq cl^d.$$

We now let $\pi_{**}$ be the $\pi_{j(x)}$ occurring most in (5.23), where $x$ varies over $\mathcal{E}'$, and define $\mathcal{E} \subseteq \mathcal{E}'$ as the subset of those $x$ in $\mathcal{E}'$ for which $\pi_{j(x)} = \pi_{**}$. With this choice, (5.14) holds by (5.24). Moreover, since (5.12) and (5.13) both hold for $\mathcal{E}'$ and since $|\mathcal{E}| \geq \frac{1}{d+1}|\mathcal{E}'|$, the same identities hold for $\mathcal{E}$ as well (with a different constant). Hence, the proof of Lemma 5.3 is complete under (5.19).

On the other hand, let us now assume (5.19) does not hold. That is, we suppose that

$$(5.25) \qquad |\mathcal{A}| > c_4|\mathcal{C}_l|.$$

We then claim that, for $N \geq c(\gamma, \gamma')$,

$$(5.26) \qquad |\{x \in \mathcal{A} : |C_x(l) \cap I| > c_4 l^{d+1}\}| \leq c_4|\mathcal{A}|.$$

Indeed, we would otherwise have

$$|I| \quad \geq \quad |\{x \in \mathcal{A} : |C_x(l) \cap I| > c_4 l^{d+1}\}|c_4 l^{d+1} \overset{\text{(if (5.26) false)}}{>} c_4^2|\mathcal{A}|l^{d+1}$$

$$\overset{(5.25)}{>} c_4^3 l^{d+1}|\mathcal{C}_l| \overset{(5.19)}{>} \frac{4}{5}l^{d+1}|\mathcal{C}_l| = \frac{4}{5}\left(\frac{l^{d+1}|\mathcal{C}_l|}{L^{d+1}}\right)L^{d+1},$$

contradicting the choice of $I$ for $N \geq c(\gamma, \gamma')$, because $\frac{l^{d+1}|\mathcal{C}_l|}{L^{d+1}}$ only depends on $N, \gamma, \gamma'$ and tends to 1 as $N \to \infty$. It follows that for $N \geq c(\gamma, \gamma')$,

$$(5.27) \quad \begin{aligned} c_4|\mathcal{C}_l| &\overset{(5.25)}{\leq} |\mathcal{A}| \overset{(5.26)}{\leq} \frac{1}{1-c_4}|\{x \in \mathcal{A} : |C_x(l) \cap I| \leq c_4 l^{d+1}\}| \\ &\overset{(5.16)}{=} \frac{1}{1-c_4}\left|\left\{x \in \mathcal{C}_l : \frac{1}{8}l^{d+1} < |C_x(l) \cap I| \leq c_4 l^{d+1}\right\}\right|. \end{aligned}$$

Defining $\mathcal{E}' = \{x \in \mathcal{C}_l : \frac{1}{8}l^{d+1} < |C_x(l) \cap I| \leq c_4 l^{d+1}\}$, we apply the isoperimetric inequality (5.11) of Lemma 5.2 with $C_x(l)$ in place of $[0, M]^{d+1}$ and $C_x(l) \cap I$ in place of $A$ for every $x \in \mathcal{E}'$ and thus obtain a projection $\pi_{j(x)}$ satisfying (5.24), as in the previous case. We then define $\mathcal{E} \subseteq \mathcal{E}'$ as the subset containing only those $x \in \mathcal{E}'$ for which $\pi_{j(x)}$ in (5.24) is equal to the most frequently occurring $\pi_{**}$. As a consequence, (5.14) holds. Moreover, (5.12) is clear by definition of $\mathcal{E}$. And finally, we have by (5.27), $|\mathcal{E}| \geq \frac{1}{d+1}|\mathcal{E}'| \geq c|\mathcal{C}_l| \geq c'(\frac{L}{l})^{d+1}$, which yields (5.13) by (5.9). This completes the proof of Lemma 5.3.  $\square$

The last geometric lemma in this section is essentially a modification of Lemma 5.3. It provides a similar result for $B(\alpha)$, $0 < d\alpha < 1$ instead of $C(L)$. The idea of the proof, illustrated in Figure 4, is to "pile up" approximately $N^{1-d\alpha}$ copies of $B(\alpha)$ into the $\mathbb{Z}$-direction of $E$ and to then apply the same arguments with the isoperimetric inequality (5.11) as in the proof of Lemma 5.3 to the resulting set intersected with $B_\infty(x, [N/4])$.

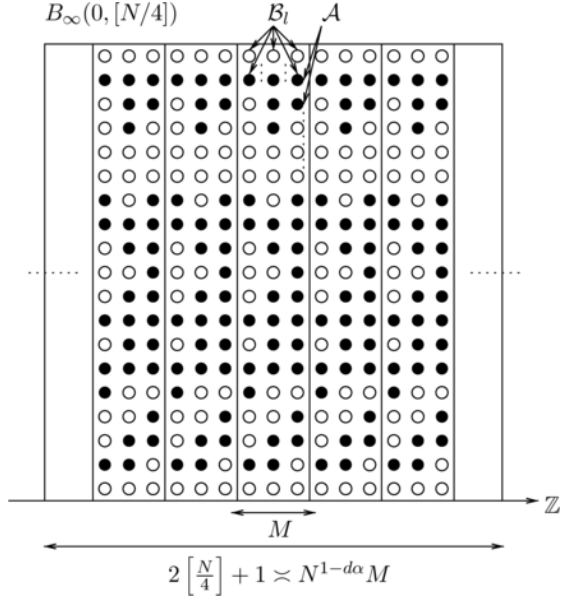

FIG. 4. *An illustration of the set $\mathcal{A}'$ of copies of $\mathcal{A} \subseteq \mathcal{B}_l$ piled up in the (horizontal) $\mathbb{Z}$-direction [cf. (5.35)], used in the proof of Lemma 5.4. The circles are the points on the subgrid $\mathcal{H}_l$ in (5.31), and the filled circles are the points contained in the set $\mathcal{A}'$. Each copy of $\mathcal{B}_l$ has thickness $M$, defined in (5.34), so that the larger box $B_\infty(0, [N/4])$ contains roughly $N^{1-d\alpha}$ copies of $\mathcal{B}_l$.*

LEMMA 5.4 ($d \geq 1$, $0 < \gamma < d\alpha < 1$, $l = [N^\gamma]$, $N \geq 1$). *For all $N \geq c(\alpha, \gamma)$, whenever $K \subseteq B(\alpha)$ [cf. (1.7)] $\frac{1}{3}$-disconnects $B(\alpha)$, there exists a set $\mathcal{E} \subseteq B(\alpha)^{(-l)}$ and projections $\pi_*$ and $\pi_{**} \in \{\pi_1, \ldots, \pi_{d+1}\}$ such that*

$$(5.28) \qquad \pi_*(\mathcal{E}) \subseteq \pi_*(B(\alpha) \cap \pi_E([-[N/4], [N/4]]^{d+1} \cap l\mathbb{Z}^{d+1})),$$

$$(5.29) \qquad |\pi_*(\mathcal{E})| \geq c'\left(\frac{N}{l}\right)^d N^{d\alpha-1},$$

$$(5.30) \qquad \text{for all } x \in \mathcal{E} : |\pi_{**}(K \cap C_x(l))| \geq c''l^d.$$

PROOF. The proof is very similar to the one of Lemma 5.3. We choose a set $I \subseteq B(\alpha)$ such that $\frac{1}{3}|B(\alpha)| \leq |I| \leq \frac{2}{3}|B(\alpha)|$ and $\partial_{B(\alpha)}(I) \subseteq K$. We then introduce the subgrids of side-length $l$ of $[-[N/4], [N/4]]^d \times \mathbb{Z}$ and of

$B(\alpha)^{(-l)}$ as [cf. (2.3)]

$$
(5.31) \qquad
\begin{aligned}
\mathcal{H}_l &= \pi_E(([-[N/4],[N/4]]^d \times \mathbb{Z}) \cap l\mathbb{Z}^{d+1}) \quad \text{and} \\
\mathcal{B}_l &= B(\alpha)^{(-l)} \cap \mathcal{H}_l,
\end{aligned}
$$

and set

$$
\mathcal{A} = \{x \in \mathcal{B}_l : |C_x(l) \cap I| > \tfrac{1}{6}l^{d+1}\}.
$$

Since the disjoint union $\bigcup_{x \in \mathcal{B}_l} C_x(l)$ contains all but at most $cN^d l$ points of $B(\alpha)$, we then have

$$
\begin{aligned}
\tfrac{1}{3}|B(\alpha)| \le |I| &\le \tfrac{1}{6}l^{d+1}|\mathcal{B}_l \setminus \mathcal{A}| + l^{d+1}|\mathcal{A}| + cN^d l \\
&\le \tfrac{1}{6}|B(\alpha)| + l^{d+1}|\mathcal{A}| + cN^d l,
\end{aligned}
$$

hence

$$
\left(\frac{1}{6} - cN^{\gamma - d\alpha}\right)\frac{|B(\alpha)|}{l^{d+1}} \le |\mathcal{A}|,
$$

and thus for $N \ge c(\alpha, \gamma)$,

$$
(5.32) \qquad c|\mathcal{B}_l| \le \frac{c|B(\alpha)|}{l^{d+1}} \le |\mathcal{A}|.
$$

Suppose now that in addition

$$
(5.33) \qquad |\mathcal{A}| \le c_5|\mathcal{B}_l| \quad \text{with } c_5^3 = \tfrac{1}{2}(1 + \tfrac{3}{4}).
$$

Then we define the set $\mathcal{A}' \subseteq \mathcal{H}_l$ by "piling up" adjoining copies of the set $\mathcal{B}_l \supseteq \mathcal{A}$ into the $\mathbb{Z}$-direction. That is, we introduce the "thickness" $M$ of $\mathcal{B}_l$,

$$
(5.34) \qquad M = \sup_{(u,v),(u',v') \in \mathcal{B}_l} |v - v'| = 2\left[\frac{[(1/4)N^{d\alpha}] - l}{l}\right]l,
$$

and define

$$
(5.35) \qquad \mathcal{A}' = \bigcup_{n \in \mathbb{Z}}(n(M + l)e_{d+1} + \mathcal{A}) \subseteq \bigcup_{n \in \mathbb{Z}}(n(M + l)e_{d+1} + \mathcal{B}_l) = \mathcal{H}_l;
$$

cf. (5.31), Figure 4. Observe that $B_\infty(0, [N/4]) \cap \mathcal{A}'$ contains no less than $cN^{1-d\alpha}$ and no more than $c'N^{1-d\alpha}$ copies of $\mathcal{A}$. With (5.32) and (5.33) it follows that for $N \ge c(\alpha, \gamma)$,

$$
c'\left(\frac{N}{l}\right)^{d+1} \le |B_\infty(0, [N/4]) \cap \mathcal{A}'| \le (1 - c')\left(\frac{N}{l}\right)^{d+1}.
$$

For $N \ge c(\gamma')$, an application of the isoperimetric inequality (5.11) of Lemma 5.2 on the subgrid $\mathcal{H}_l$ defined in (5.31), with $B_\infty(0, [N/4]) \cap \mathcal{H}_l$ in place of $[0, M]^{d+1}$ and $B_\infty(0, [N/4]) \cap \mathcal{A}'$ in place of $A$, hence yields an $i \in \{1, \ldots, d+1\}$ such that

$$(5.36) \qquad |\pi_i(\partial_{\mathcal{H}_l, i}(\mathcal{A}'))| \geq c \left(\frac{N}{l}\right)^d.$$

If $i \neq d+1$, then the set on the left-hand side of (5.36) is contained in the at most $cN^{1-d\alpha}$ translated copies of the set $\pi_i(\partial_{\mathcal{B}_l, i}(\mathcal{A}))$ intersecting $B_\infty(0, [N/4])$ [see (5.35) and Figure 4]. We then deduce from (5.36) that

$$(5.37) \qquad |\pi_i(\partial_{\mathcal{B}_l, i}(\mathcal{A}))| \geq c \left(\frac{N}{l}\right)^d N^{d\alpha - 1}.$$

If $i = d+1$, in (5.36), then we claim that

$$(5.38) \qquad \pi_{d+1}(\partial_{\mathcal{B}_l, d+1}(\mathcal{A})) \supseteq \pi_{d+1}(\partial_{\mathcal{H}_l, d+1}(\mathcal{A}')).$$

Indeed, suppose some $u \in \mathbb{T}_N^d$ does not belong to the left-hand side. Then the fiber $\{x \in \mathcal{B}_l : \pi_{d+1}(x) = u\}$ must either be disjoint from $\mathcal{A}$ or be a subset of $\mathcal{A}$. Our construction of $\mathcal{A}'$ in (5.35) implies that the set $\{x \in \mathcal{H}_l : \pi_{d+1}(x) = u\}$ is then either disjoint from $\mathcal{A}'$ or a subset of $\mathcal{A}'$, as in the first and second horizontal lines of Figure 4 [note that the translated copies of $\mathcal{B}_l$ in (5.35) adjoin each other on the subgrid $\mathcal{H}_l$]. But this precisely means that $u$ is not included in the right-hand side of (5.38). In particular, by (5.36) and (5.38), (5.37) holds also with $i = d+1$ (even without the $N^{d\alpha - 1}$ on the right-hand side). Using (5.37), we can perform the same construction as in the proof of Lemma 5.3 below (5.20) in order to obtain the desired set $\mathcal{E}$.

If, on the other hand, (5.33) does not hold, that is, if

$$|\mathcal{A}| > c_5 |\mathcal{B}_l|,$$

then the existence of the required set $\mathcal{E}$ follows from the argument below (5.25), where (5.29) can be deduced from $|\mathcal{E}| \geq c|\mathcal{B}_l| \geq c'(N/l)^{d+1} N^{d\alpha - 1}$ by applying the estimate (5.9) to $[cN^{1-d\alpha}]$ copies of $\mathcal{E}$ piled-up in a box. $\quad \square$

## 6. The large deviations estimate.

Our task in this last section is to derive the following form of the large deviations estimate (1.10):

THEOREM 6.1 ($d \geq 3$). *The estimate (1.10) holds with (cf. Figure 2)*

$$(6.1) \quad f(\alpha, \beta) = \begin{cases} d - 1 - \dfrac{d\alpha}{d-1}, & \text{on } (0, 1/d) \times \left(0, d-1-\dfrac{d\alpha}{d-1}\right), \\ d - 1 - \dfrac{1}{d-1}, & \text{on } [1/d, \infty) \times \left(0, d-1-\dfrac{1}{d-1}\right), \\ ((d-1)^2 - 1)(d-1-\beta), & \\ & \text{on } [1/d, \infty) \times \left[d-1-\dfrac{1}{d-1}, d-1\right), \\ 0, & \text{otherwise.} \end{cases}$$

Before we begin with the proof of Theorem 6.1, we examine its implications. With the function $f$ in (6.1), the lower bound exponents $d(1 - \alpha - \varphi(\alpha))$ [in (1.4)] and $\zeta$ [in (1.13)] are related via (1.14), as will be checked in Corollary 6.3. We therefore have to justify the expression $\vee d(1 - 2\alpha)\mathbf{1}_{\{\alpha < 1/d\}}$ on the right-hand side of (1.14). This is the aim of the next proposition.

PROPOSITION 6.2 $(d \geq 2, \ 0 < \alpha < \frac{1}{d})$.   *For some constant $c_6 > 0$,*

$$(6.2) \qquad P_0^{N^{-d\alpha}}[\exp\{c_6 N^{d(1-2\alpha)}\} \leq T_N^{\mathrm{disc}}] \xrightarrow{N \to \infty} 1.$$

PROOF.   The idea is that, by our previous geometric estimates, any trajectory disconnecting $E$ must contain at least $cN^d$ points in a box of the form $x + B_\infty(0, [N/4])$, $x \in E$. Hence, there must be two visited points within distance $N$ from each other, such that the random walk $X$ spends $[cN^d]$ time units between the visits to the two points. The probability of this event can be bounded from above by standard large deviations estimates.

In detail: Lemma 4.1, applied with $B(\alpha) = B_\infty(0, [N/4])$ (i.e., with $\alpha \geq \frac{1}{d}$), shows that, for $t \geq 0$, $N \geq c$, the event $\{X([0, [t]])$ disconnects $E\}$ is contained in the event

$$(6.3) \qquad \bigcup_{\substack{x \in E \\ |x_{d+1}| \leq [t]+N}} \{X([0, [t]]) \ \tfrac{1}{3}\text{-disconnects } x + B_\infty(0, [N/4])\}.$$

We now choose a set $I \subseteq x + B_\infty(0, [N/4])$ corresponding to $\frac{1}{3}$-disconnection of $x + B_\infty(0, [N/4])$ by $X([0, [t]])$ [cf. (1.6)]. By the isoperimetric inequality (5.11) of Lemma 5.2, applied with $x + B_\infty(0, [N/4])$ in place of $[0, M]^{d+1}$ and $I$ in place of $A$, the event (6.3) is contained in $\bigcup_{\substack{x \in E \\ |x_{d+1}| \leq [t]+N}} A_x([t])$, where, for some constant $c_7 > 0$,

$$A_x([t]) = \{|X([0, [t]]) \cap (x + B_\infty(0, [N/4]))| \geq c_7 N^d\}.$$

We therefore have

$$
\begin{aligned}
P_0^{N^{-d\alpha}}[T_N^{\mathrm{disc}} \leq t] &\leq P_0^{N^{-d\alpha}}\Bigg[\bigcup_{\substack{x \in E \\ |x_{d+1}| \leq [t]+N}} A_x([t])\Bigg] \\
(6.4) \qquad &\leq cN^d(t+N)\sup_{x \in E} P_0^{N^{-d\alpha}}[A_x([t])].
\end{aligned}
$$

By the strong Markov property applied at $H_{x+B_\infty(0,[N/4])}^X$, the entrance time of $x + B_\infty(0, [N/4])$, and using translation invariance of $X$, we obtain

$$\sup_{x \in E} P_0^{N^{-d\alpha}}[A_x([t])]$$

$$\leq \sup_{x \in E} P_0^{N^{-d\alpha}}[\theta_{H^X_{(x+B_\infty(0,[N/4]))}}^{-1} A_x([t])]$$

$$\overset{\text{(Markov, transl. inv.)}}{\leq} \sup_{x\,:\,A_x([t]) \ni 0} P_0^{N^{-d\alpha}}[A_x([t])]$$

$$\leq P_0^{N^{-d\alpha}}[\text{for some } n \geq c_7 N^d : \pi_{\mathbb{Z}}(X_n) \leq N].$$

Inserting this last inequality into (6.4) and using that $N \leq \frac{N^{-d\alpha}n}{2(d+1)}$ for $n \geq c_7 N^d$, $N \geq c(\alpha)$ (because $d - d\alpha > d - 1 \geq 1$), we deduce that, for $N \geq c(\alpha)$,

$$
\begin{aligned}
(6.5) \quad & P_0^{N^{-d\alpha}}[T_N^{\text{disc}} \leq t] \\
& \quad \leq cN^d(t+N) P_0^{N^{-d\alpha}}[\text{for some } n \geq c_7 N^d : \pi_{\mathbb{Z}}(X_n) \leq N] \\
& \quad \leq cN^d(t+N) \sum_{n \geq c_7 N^d} P_0^{N^{-d\alpha}}\left[\pi_{\mathbb{Z}}(X_n) \leq \frac{N^{-d\alpha}n}{2(d+1)}\right] \\
& \quad = cN^d(t+N) \sum_{n \geq c_7 N^d} P_0^{N^{-d\alpha}}\left[\pi_{\mathbb{Z}}(X_n) - \frac{N^{-d\alpha}n}{d+1} < -\frac{N^{-d\alpha}n}{2(d+1)}\right].
\end{aligned}
$$

Since $(\pi_{\mathbb{Z}}(X_n) - \frac{N^{-d\alpha}n}{d+1})_{n \geq 0}$ is a $P_0^{N^{-d\alpha}}$-martingale with steps bounded by $c$, Azuma's inequality (cf. [2], page 85) implies that

$$P_0^{N^{-d\alpha}}\left[\pi_{\mathbb{Z}}(X_n) - \frac{N^{-d\alpha}n}{d+1} < -\frac{N^{-d\alpha}n}{2(d+1)}\right] \leq e^{-cN^{-2d\alpha}n}.$$

Applying this estimate to (6.5) with $t_N = \exp\{c_6 N^{d-2d\alpha}\}$ we see that for $N \geq c(\alpha)$,

$$P_0^{N^{-d\alpha}}[T_N^{\text{disc}} \leq \exp\{c_6 N^{d-2d\alpha}\}] \leq cN^{d+2d\alpha} \exp\{c_6 N^{d-2d\alpha} - c' N^{d-2d\alpha}\}.$$

Choosing the constant $c_6 > 0$ sufficiently small, this yields (6.2) (recall that $d\alpha < 1 \leq \frac{d}{2}$). □

We can now check that Theorem 6.1 does have the desired implications on the lower bounds on $T_N^{\text{disc}}$.

COROLLARY 6.3 ($d \geq 3$, $\alpha > 0$, $\varepsilon > 0$).   With $\varphi$ defined in (1.5), one has

$$(6.6) \qquad \text{for } \alpha > 1, \qquad\qquad P_0^{N^{-d\alpha}}[N^{2d-\varepsilon} \leq T_N^{\text{disc}}] \overset{N \to \infty}{\longrightarrow} 1,$$

$$(6.7) \qquad \text{for } \alpha < 1, \qquad P_0^{N^{-d\alpha}}[\exp\{N^{d(1-\alpha-\varphi(\alpha))-\varepsilon}\} \leq T_N^{\text{disc}}] \overset{N \to \infty}{\longrightarrow} 1.$$

PROOF.  Since the function $f$ of (6.1) satisfies $f(\alpha, \beta) > 0$ for $(\alpha, \beta) \in (1, \infty) \times (0, d-1)$, (6.6) follows immediately from Theorem 6.1 and (1.11).

By (1.12) and (6.2), (6.7) holds with $\varphi$ defined for $\alpha \in (0,1)$ by (1.14). Let us check that the expression for $\varphi$ in (1.14) agrees with (1.5). We first treat the case $\alpha \in [\frac{1}{d}, 1)$, for which $f(\alpha, \cdot)$ is illustrated on the right-hand side of Figure 2, below Theorem 1.2. We have $d\alpha \geq 1$, $f(\alpha, \beta) = 0$ for $\beta \geq d - 1$ and the maximum of $g_\alpha$ [cf. (1.13)] on $(0, d-1)$ is attained at (see Figure 2)

$$\bar\beta = d - 1 - \frac{d - d\alpha}{(d-1)^2} \in \left[ d - 1 - \frac{1}{d-1}, d-1 \right) \cap (d\alpha - 1, d-1).$$

Hence, for $\alpha \in [\frac{1}{d}, 1)$ [cf. (1.13)],

$$\zeta = \sup_{\beta > 0} g_\alpha(\beta) = g_\alpha(\bar\beta) = d\left( 1 - \alpha - \frac{1-\alpha}{(d-1)^2} \right),$$

and therefore $\varphi(\alpha) = \frac{1-\alpha}{(d-1)^2}$, as required.

Turning to the case $\alpha \in (0, \frac{1}{d})$, we refer to the left-hand side of Figure 2 for an illustration of $f$. We now have $d\alpha - 1 < 0$, $f(\alpha, \beta) = 0$ for $\beta > d - 1 - \frac{d\alpha}{d-1}$ and hence

$$\zeta = \sup_{\beta > 0} g_\alpha(\beta)$$

$$= \sup_{\beta \in (0, d-1-d\alpha/(d-1))} \left( \beta \wedge d - 1 - \frac{d\alpha}{d-1} \right) = d\left( 1 - \frac{1}{d} - \frac{\alpha}{d-1} \right).$$

Therefore [cf. (1.14)], for $\alpha \in (0, \frac{1}{d})$,

$$(6.8) \qquad \varphi(\alpha) = 1 - \alpha - \left( \left( 1 - \frac{1}{d} - \frac{\alpha}{d-1} \right) \vee (1 - 2\alpha) \right).$$

This expression is immediately seen to coincide with (1.5) for $\alpha \in (0, \frac{1}{d})$ near 0 and $\frac{1}{d}$, and $\alpha_*$ is precisely the value for which $1 - \frac{1}{d} - \frac{\alpha_*}{d-1} = 1 - 2\alpha_*$, so that (6.8), and hence (1.14), agrees with (1.5).  □

Thanks to Corollary 6.3, the lower bounds on $T_N^{\mathrm{disc}}$ of Theorem 1.1 will be established once we show Theorem 6.1. Let us give a rough outline of the strategy of the proof. In the previous section, we have shown that if $K = X([0, \mathcal{D}_{[N^\beta]}])$ $\frac{1}{3}$-disconnects $B(\alpha)$, then there must be a wealth of subcubes of $B(\alpha)$ such that $X([0, \mathcal{D}_{[N^\beta]}])$ contains a surface of points in every subcube (see Lemmas 5.3 and 5.4 for the precise statements and Figure 3 for an illustration). The crucial upper bound on the probability of an event of this form is obtained in Lemma 6.5, using Khaśminskii's lemma to obtain an exponential tail estimate on the number of points visited by $X$ during a suitably defined excursion. This upper bound is then applied in order to find the needed large deviations estimate of the form (1.10). We begin by collecting the required estimates involving the Green function [cf. (2.25)].

LEMMA 6.4 ($d \geq 2$, $N, a \geq 1$, $100 \leq a \leq 4N$, $A \subseteq B \subseteq S_a$).

$$(6.9) \qquad P_x^0[H_A^X < H_{B^c}^X] \leq \frac{\sum_{y \in A} g^B(x, y)}{\inf_{y \in A} \sum_{y' \in A} g^B(y, y')} \qquad for \; x \in B.$$

*For any $x, x' \in S_a$, one has*

$$(6.10) \qquad g^{S_a}(x, x') \leq c(1 \vee |x - x'|_\infty)^{1-d} \exp\left\{-c' \frac{|x - x'|_\infty}{a}\right\}.$$

*If $\operatorname{diam}(A) \leq \frac{a}{100}$ [cf. (2.4)] and $A \subseteq B^{-(a/10)}$ [cf. (2.3)], then, for $x, x' \in A$,*

$$(6.11) \qquad c|x - x'|_\infty^{1-d} \leq g^B(x, x').$$

PROOF. The estimate (6.9) follows from an application of the strong Markov property at $H_A^X$. The estimate (6.10) follows from the bound on the Green function of the simple random walk on $\mathbb{Z}^{d+1}$ killed when exiting the slab $\mathbb{Z}^d \times [-[a], [a]]$ in (2.13) of Sznitman [11]. For (6.11), we note that, by assumption, $B_\infty(x, \frac{a}{10}) \subseteq B$. In particular, it follows from translation invariance that

$$(6.12) \qquad g^B(x, x') \geq g^{B_\infty(0, a/10)}(0, x - x').$$

By assumption $\frac{a}{10} \leq \frac{2N}{5}$, so the right-hand side of (6.12) can be identified with the corresponding Green function for the simple random walk on $\mathbb{Z}^{d+1}$, and (6.11) follows from the estimate of Lawler [9], page 35, Proposition 1.5.9. □

We now introduce, for sets $U, \tilde{U} \subseteq E$, the times $(\tilde{R}_n)_{n \geq 1}$ and $(\tilde{D}_n)_{n \geq 1}$ as the times of return to $U$ and departure from $\tilde{U}$ [cf. (2.24)] and denote with $\pi_*$ and $\pi_{**}$ elements of the set of projections $\{\pi_1, \dots, \pi_{d+1}\}$. The next lemma then provides a control on an event of the form [cf. (2.3), (2.8)]

$$(6.13)
\begin{aligned}
&A_{U, \tilde{U}, l, M_1, M_2} \\
&= \bigcup_{\pi_*, \pi_{**}} \bigcup_{\substack{\mathcal{E} \subseteq U^{(-l)} \\ |y - y'|_\infty \geq l \text{ for } y, y' \in \mathcal{E}, \\ |\pi_*(\mathcal{E})| \geq M_1}} \bigcap_{y \in \mathcal{E}} \{|\pi_{**}(X([0, \tilde{D}_{M_2}]) \cap C_y(l))| \geq cl^d\}.
\end{aligned}$$

Our method does not produce a useful upper bound when $d = 2$ [note that when $d = 2$, the right-hand side of (6.14) is greater than 1 for $N \geq c$]. Although it is possible to obtain a bound for $d = 2$ tending to 0 as $N \to \infty$, using estimates on the Green function in dimension 2, it does not seem to be possible to obtain an exponential decay in $N$ with this approach. Thus, the upper bound we have for $d = 2$ brings little information on the large deviations problem (1.10).

LEMMA 6.5 ($d \geq 3$, $N, l, a, M_1, M_2 \geq 1$, $100 \leq a \leq 4N$, $1 \leq l \leq \frac{a}{100}$). Let $U$, $\tilde{U} \subseteq E$ be sets such that $U \subseteq U^{(a/10)} \subseteq \tilde{U} \subseteq x_* + S_a$ [cf. (2.2), (1.9)]. Then one has the estimate

$$(6.14) \qquad \sup_{x \in E} P_x^0[A_{U, \tilde{U}, l, M_1, M_2}] \leq \exp\{c' M_2 + c' M_1 \log N - c'' M_1 a^{-1} l^{d-1}\}$$

[on the event defined in (6.13)].

PROOF. In order to abbreviate the notation, we denote the event in (6.13) by $A$ during the proof. Furthermore, by replacing $\mathcal{E}$ with a subset, we may assume that

$$(6.15) \qquad |\pi_*(\mathcal{E})| = |\mathcal{E}| = M_1.$$

Also, translation invariance allows us to set $x_* = 0$.

The first step is to note that the number of possible choices of the set $\mathcal{E}$ in the definition of $A$ is not larger than

$$|U|^{|\mathcal{E}|} \overset{(6.15)}{\leq} (cN)^{(d+1)M_1} \leq \exp\{c M_1 \log N\}.$$

Next, we note that visits made by the random walk $X$ to $C_y(l)$, $y \in \mathcal{E}$, can only occur during the time intervals $[\tilde{R}_n, \tilde{D}_n]$, $n \geq 1$ (because $\mathcal{E} \subseteq U^{(-l)}$). From these observations, we deduce that

$$
\begin{aligned}
(6.16) \qquad & \sup_{x \in E} P_x^0[A] \\
& \leq c e^{c M_1 \log N} \\
& \quad \times \sup_{x, \mathcal{E}, \pi_*, \pi_{**}} P_x^0 \left[ \sum_{n=1}^{M_2} \sum_{y \in \mathcal{E}} |\pi_{**}(X([\tilde{R}_n, \tilde{D}_n]) \cap C_y(l))| \geq c M_1 l^d \right],
\end{aligned}
$$

where the supremum is taken over all $x \in E$, and all possible sets $\mathcal{E}$ and projections $\pi_*, \pi_{**}$ entering the definition of the event $A$. By the exponential Chebyshev inequality and the strong Markov property applied inductively at $\tilde{R}_{M_2}, \tilde{R}_{M_2-1}, \ldots, \tilde{R}_1$, it follows from (6.16) that, for any $r \geq 1$, $\sup_{x \in E} P_x^0[A]$ is bounded by

$$
\begin{aligned}
(6.17) \qquad & c e^{c M_1 \log N - c r M_1 l^d} \\
& \quad \times \sup_{x, \mathcal{E}, \pi_*, \pi_{**}} E_x^0 \left[ \exp\left\{ \sum_{n=1}^{M_2} \sum_{y \in \mathcal{E}} r |\pi_{**}(X([0, \tilde{D}_1]) \cap C_y(l))| \circ \theta_{\tilde{R}_n} \right\} \right] \\
& \overset{\text{(Markov)}}{\leq} c e^{c M_1 \log N - c r M_1 l^d} \\
& \quad \times \sup_{\mathcal{E}, \pi_*, \pi_{**}} \left( \sup_{x \in U} E_x^0 \left[ \exp\left\{ \sum_{y \in \mathcal{E}} r |\pi_{**}(X([0, \tilde{D}_1]) \cap C_y(l))| \right\} \right] \right)^{M_2}.
\end{aligned}
$$

Before deriving an upper bound on this last expectation, we introduce the following notational simplification: for any point $z \in C_y(l)$, we denote its fiber in $C_y(l)$ of points of equal $\pi_{**}$-projection by $J_z$, or in other words, for $z \in C_y(l)$,

$$J_z = \{z' \in C_y(l) : \pi_{**}(z') = \pi_{**}(z)\}.$$

The collection of all fibers in the box $C_y(l)$ is denoted by

$$(6.18) \qquad F(y) = \{J_z : z \in C_y(l)\},$$

and the collection of all fibers by

$$(6.19) \qquad F = \bigcup_{y \in \mathcal{E}} F(y).$$

Using this notation, we have [cf. (2.22)]

$$(6.20) \qquad \sum_{y \in \mathcal{E}} |\pi_{**}(X([0, \tilde{D}_1]) \cap C_y(l))| = \sum_{J \in F} \mathbf{1}_{\{H_J^X < \tilde{D}_1\}}.$$

By the version of Khaśminskii's lemma of (2.46) of Dembo and Sznitman [4] (see also [8]), we see that for any $x \in U$ and $r \geq 0$,

$$(6.21) \quad E_x^0\left[\exp\left\{r \sum_{J \in F} \mathbf{1}_{\{H_J^X < \tilde{D}_1\}}\right\}\right] \leq \sum_{k \geq 0} r^k \left(\sup_{x \in U} E_x^0\left[\sum_{J \in F} \mathbf{1}_{\{H_J^X < \tilde{D}_1\}}\right]\right)^k.$$

Writing [cf. (6.18), (6.19)]

$$\sum_{J \in F} \mathbf{1}_{\{H_J^X < \tilde{D}_1\}} = \sum_{y \in \mathcal{E}} \sum_{J \in F(y)} \mathbf{1}_{\{H_J^X < \tilde{D}_1\}},$$

for any $x \in U$, the strong Markov property applied at $H_{C_y(l)}^X$ yields

$$(6.22) \quad \begin{aligned} &E_x^0\left[\sum_{J \in F} \mathbf{1}_{\{H_J^X < \tilde{D}_1\}}\right] \\ &= \sum_{y \in \mathcal{E}} E_x^0\left[H_{C_y(l)}^X < \tilde{D}_1, \sum_{J \in F(y)} (\mathbf{1}_{\{H_J^X < \tilde{D}_1\}}) \circ \theta_{H_{C_y(l)}^X}\right] \\ &\leq \sum_{y \in \mathcal{E}} P_x^0[H_{C_y(l)}^X < \tilde{D}_1] \sup_{z \in C_y(l)} E_z^0\left[\sum_{J \in F(y)} \mathbf{1}_{\{H_J^X < \tilde{D}_1\}}\right]. \end{aligned}$$

To bound the right-hand side of (6.22), we note that, for any $z \in C_y(l)$ and $k \in \{0, \ldots, l-1\}$, at most $c(1 \vee k)^{d-1}$ of the fibers $J \in F(y)$ are at $|\cdot|_\infty$-distance $1 \vee k$ from $J_z$ and thus deduce that, for any $z \in C_y(l)$,

$$E_z^0\left[\sum_{J \in F(y)} \mathbf{1}_{\{H_J^X < \tilde{D}_1\}}\right]$$

(6.23)

$$\leq c \sum_{k=0}^{l-1} (1 \vee k)^{d-1} \sup_{z' \, : \, |\pi_{**}(z-z')|_\infty = k} P_z^0[H_{J_{z'}}^X < \tilde{D}_1].$$

For this last probability, we use the estimate (6.9), applied with $A = J_{z'}$, $B = \tilde{U}$ and $x = z$. With the help of (6.10) and the assumption that $\tilde{U} \subseteq S_a$, the numerator of the right-hand side of (6.9) can then be bounded from above by $clk^{1-d}$, while the denominator is trivially bounded from below by 1. We thus obtain

$$\sup_{z' \, : \, |\pi_{**}(z-z')| = k} P_z^0[H_{J_{z'}}^X < \tilde{D}_1] \leq clk^{1-d}.$$

With (6.23), this yields

$$E_z^0\left[\sum_{J \in F(y)} \mathbf{1}_{\{H_J^X < \tilde{D}_1\}}\right] \leq cl^2 \qquad \text{for any } z \in C_y(l).$$

Coming back to (6.22), we obtain

(6.24) $\quad E_x^0\left[\sum_{J \in F} \mathbf{1}_{\{H_J^X < \tilde{D}_1\}}\right] \leq cl^2 \sum_{y \in \mathcal{E}} P_x^0[H_{C_y(l)}^X < \tilde{D}_1] \qquad \text{for any } x \in U.$

For this last sum, we proceed as before: Note that, by (6.15), the sum can be regarded as a sum over the set $\pi_*(\mathcal{E})$, which is a subset of the $d$-dimensional lattice $\pi_*(E)$. Since moreover $|y - y'|_\infty \geq l$ for all $y, y' \in \mathcal{E}$, there are at most $c(1 \vee k)^{d-1}$ points in $\pi_*(\mathcal{E})$ of $|\cdot|_\infty$-distance between $kl$ and $(k+1)l$ from $\pi_*(x)$. We therefore have, for any $x \in U$,

$$\sum_{y \in \mathcal{E}} P_x^0[H_{C_y(l)}^X < \tilde{D}_1]$$

(6.25)

$$\leq c \sum_{k=0}^{\infty} (1 \vee k)^{d-1} \sup_{y \in \mathcal{E} \, : \, |\pi_*(y-x)| \geq kl} P_x^0[H_{C_y(l)}^X < \tilde{D}_1].$$

In order to bound this last probability, we again use the estimate (6.9), this time with $A = C_y(l)$ and $B = \tilde{U}$. By (6.10), our assumption that $\tilde{U} \subseteq S_a$ then allows us to bound the numerator of the right-hand side of (6.9) from above by $cl^{d+1}(1 \vee lk)^{1-d}e^{-c'lk/a}$, while our assumptions $l \leq \frac{a}{100}$ and $C_y(l) \subseteq U \subseteq \tilde{U}^{(-a/10)}$ allow us to use (6.11) and find the lower bound of $cl^2$ on the denominator. We thus have

$$\sup_{y \in \mathcal{E} \, : \, |\pi_*(y-x)| \geq lk} P_x^0[H_{C_{x'}(l)}^X < \tilde{D}_1] \leq c(1 \vee k)^{1-d}e^{-c'(lk/a)}.$$

With (6.25), this yields

$$\sum_{y \in \mathcal{E}} P_x^0[H_{C_y(l)}^X < \tilde{D}_1] \leq c\frac{a}{l} \qquad \text{for any } x \in U,$$

which we insert into (6.24) to obtain

$$\sup_{x \in U} E_x^0 \left[ \sum_{J \in F} \mathbf{1}_{\{H_J^X < \tilde{D}_1\}} \right] \le c_8 a l.$$

Choosing $r = \frac{1}{2c_8 a l}$ in (6.21), we infer that

$$E_x^0 \left[ \exp \left\{ \frac{1}{2c_8 a l} \sum_{J \in F} \mathbf{1}_{\{H_J^X < \tilde{D}_1\}} \right\} \right] \le 2 \qquad \text{for any } x \in U.$$

Coming back to (6.17) with $r$ as above and remembering (6.20), we deduce (6.14) and thus complete the proof of Lemma 6.5. $\square$

The remaining part of the proof of Theorem 6.1 is essentially an application of Lemma 6.5 together with the geometric Lemmas 5.1–5.4 showing that the event on the left-hand side of (1.10) is contained in a union of events of the form (6.13). For $\alpha < \frac{1}{d}$, all that remains to be done is to combine Lemma 5.5 with Lemma 6.5. For $\alpha \ge \frac{1}{d}$, that is, for $B(\alpha) = B_\infty(0, [N/4])$, we additionally use an upper bound on the probability that the random walk $X$ makes a certain number of excursions between $C_{x_*}(L)$ and $(C_{x_*}(L)^{(L)})^c$ until time $\mathcal{D}_{[N^\beta]}$ for $x_* \in B_\infty(0, [N/4])$ and $L$ as in (5.1) (cf. Lemma 6.6) before we apply the geometric Lemmas 5.1 and 5.3 and the estimate (6.14) with $U = C_{x_*}(L)$.

PROOF OF THEOREM 6.1—CASE $\alpha < \frac{1}{d}$. In this case, we have to show (1.10) with $f$ illustrated on the left-hand side of Figure 2 (below Theorem 1.2) and $[N^{d\alpha \wedge 1}] = [N^{d\alpha}]$. Lemma 5.4 implies that, for $l$ as in (5.1) and the event $A_{\cdot,\cdot,\cdot,\cdot}$ defined in (6.13),

(6.26) $\quad \{U_{B(\alpha)} \le \mathcal{D}_{[N^\beta]}\} \subseteq A_{S_{2[N^{d\alpha}]}, S_{4[N^{d\alpha}]}, l, cN^{d-1+d\alpha l - d}, [N^\beta]} \overset{\text{(def.)}}{=} A_N'.$

Lemma 6.5, applied with $a = 4[N^{d\alpha}]$ and $x_* = 0$, yields

(6.27) $\quad \sup_{x \in S_{2N^{d\alpha}}} P_x^0[A_N'] \le \exp\{cN^\beta + cN^{d-1+d\alpha - d\gamma} \log N - c'N^{d-1-\gamma}\}.$

In view of (5.1), we have $0 < \gamma < d\alpha$, and provided $d-1+d\alpha-d\gamma < d-1-\gamma$ and $\beta < d-1-\gamma$, that is, if

(6.28) $\qquad \frac{d\alpha}{d-1} < \gamma < d\alpha, \qquad \beta < d-1-\gamma,$

then (6.26) and (6.27) together show that

(6.29) $\qquad \sup_{x \in S_{2[N^{d\alpha}]}} P_x^0[U_{B(\alpha)} \le \mathcal{D}_{[N^\beta]}] \le \exp\{-cN^{d-1-\gamma}\}.$

For $\beta \in (0, d - 1 - \frac{d\alpha}{d-1})$, $d \geq 3$, the constraints (6.28) are satisfied by $\gamma_0 = \frac{d\alpha}{d-1} + \varepsilon_0(d, \beta)$ for $\varepsilon_0(d, \beta) > 0$ sufficiently small. Moreover, $d - 1 - \gamma_0 = d - 1 - \frac{d\alpha}{d-1} - \varepsilon_0(d, \beta) \overset{(6.1)}{=} f(\alpha, \beta) - \varepsilon_0(d, \beta)$. Since we can make $\varepsilon_0(d, \beta) > 0$ arbitrarily small, (6.29) thus shows (1.10) for the case $\alpha \in (0, \frac{1}{d})$. This completes the proof of Theorem 6.1 in the case $\alpha < \frac{1}{d}$. $\square$

PROOF OF THEOREM 6.1—CASE $\alpha \geq \frac{1}{d}$. Recall that in this case we have to find an estimate of the form (1.10) with the function $f$ illustrated on the right-hand side of Figure 2 (below Theorem 1.2) and with $B(\alpha) = B_\infty(0, [N/4])$ [cf. (1.7)]. In order to apply Lemma 6.5 with $U = C_{x_*}(L) \subseteq B(\alpha)$, $L$ as in (5.1), we consider

(6.30)
$$\tilde{R}_n^{x_*}, \tilde{D}_n^{x_*}, n \geq 1, \text{ the successive returns to } C_{x_*}(L)$$

and departures from $C_{x_*}(L)^{(L)}$ [cf. (2.2), (2.24)].

The following lemma, in the spirit of Dembo and Sznitman [4], provides an estimate on the number of excursions between $C_{x_*}(L)$ and $(C_{x_*}(L)^{(L)})^c$ occurring during the $[N^\beta]$ excursions under consideration in (1.10).

LEMMA 6.6 ($d \geq 2$, $\alpha \geq \frac{1}{d}$, $\beta > 0$, $\gamma' \in (0, 1)$, $L = [N^{\gamma'}]$, $m, N \geq 1$). For $x \in S_{2N}$, $x_* \in B_\infty(0, [N/4])$ and $\tilde{R}_m^{x_*}$ defined in (6.30),

(6.31)
$$P_x^0[\tilde{R}_m^{x_*} \leq \mathcal{D}_{[N^\beta]}] \leq c \exp\{cN^{1-d}L^{d-1}N^\beta - c'm\}.$$

PROOF. We follow the proof of Lemma 2.3 by Dembo and Sznitman [4]. Since $C_{x_*}(L) \subseteq S_{2N}$, visits made by $X$ to $C_{x_*}(L)$ can only occur during the time intervals $[\mathcal{R}_i, \mathcal{D}_i]$, $i \geq 1$; cf. above (1.10). Let us denote the number of excursions between $C_{x_*}(L)$ and $(C_{x_*}(L)^{(L)})^c$ performed by $X$ during $[\mathcal{R}_i, \mathcal{D}_i]$ by $\mathcal{N}_i$, that is,

$$\mathcal{N}_i = |\{n \geq 1 : \mathcal{R}_i \leq \tilde{R}_n^{x_*} \leq \mathcal{D}_i\}|, \qquad i \geq 1.$$

Note that one then has $\mathcal{N}_i = \mathcal{N}_1 \circ \theta_{\mathcal{R}_i}$, $i \geq 1$. For any $\lambda > 0$, $x \in S_{2N}$, $x_* \in B_\infty(0, [N/4])$, we apply the strong Markov property at $\mathcal{R}_2$ and deduce that

$$P_x^0[\tilde{R}_m^{x_*} \leq \mathcal{D}_{[N^\beta]}]$$
$$\leq P_x^0\left[\left\{\mathcal{N}_1 \geq \left[\frac{m}{2}\right]\right\} \cup \theta_{\mathcal{R}_2}^{-1}\left\{\mathcal{N}_1 + \cdots + \mathcal{N}_{[N^\beta]-1} \geq \left[\frac{m}{2}\right]\right\}\right]$$
$$\leq P_x^0\left[\mathcal{N}_1 \geq \left[\frac{m}{2}\right]\right] + \sup_{x \in S_{2N} : |x_{d+1}| = 2N} P_x^0\left[\sum_{i=1}^{[N^\beta]-1} \mathcal{N}_i \geq \left[\frac{m}{2}\right]\right].$$

With the strong Markov property applied inductively at $\mathcal{R}_{[N^\beta]-1}, \mathcal{R}_{[N^\beta]-2}, \ldots, \mathcal{R}_1$ to the second term on the right-hand side, one infers that

$$(6.32) \quad \begin{aligned} &P_x^0[\tilde{R}_m^{x_*} \leq \mathcal{D}_{[N^\beta]}] \\ &\leq e^{-\lambda[m/2]}\left(E_x^0[e^{\lambda\mathcal{N}_1}] + \sup_{x \in S_{2N}:\, |x_{d+1}|=2N} E_x^0[e^{\lambda\mathcal{N}_1}]^{([N^\beta]-1)}\right). \end{aligned}$$

For any $x \in S_{2N}$,

$$(6.33) \quad E_x^0[e^{\lambda\mathcal{N}_1}] = 1 + (e^\lambda - 1)\sum_{n \geq 0} e^{\lambda n} P_x^0[\mathcal{N}_1 > n].$$

Applying the strong Markov property and the invariance principle as in [4], (2.16) and below, we find that

$$(6.34) \quad P_x^0[\mathcal{N}_1 > n] \leq (1-c)^n P_x^0[\mathcal{N}_1 > 0].$$

Choosing $\lambda > 0$ such that $e^\lambda(1-c) < 1$ with $c$ as in (6.34), and coming back to (6.33), we see that for any $x \in S_{2N}$,

$$(6.35) \quad E_x^0[e^{\lambda\mathcal{N}_1}] \leq 1 + c(\lambda) P_x^0[\mathcal{N}_1 > 0].$$

If we consider $|x_{d+1}| = 2N$, then we can apply the estimate (6.9) to $P_x^0[\mathcal{N}_1 > 0] = P_x^0[H_{C_{x_*}(L)}^X < \mathcal{D}_1]$ with $A = C_{x_*}(L)$, $B = S_{4N}$, $a = 4N$ and then use the Green function estimates (6.10) for the numerator and (6.11) for the denominator of the right-hand side of (6.9), to obtain, for $N \geq c(\gamma')$,

$$P_x^0[\mathcal{N}_1 > 0] \leq cL^{d-1}N^{1-d}.$$

With (6.35), this yields, for any $x \in S_{2N}$ with $|x_{d+1}| = 2N$,

$$(6.36) \quad E_x^0[e^{\lambda\mathcal{N}_1}] \leq 1 + c(\lambda)L^{d-1}N^{1-d}.$$

By (6.35), the first expectation on the right-hand side of (6.32) is bounded by a constant and with (6.36), the second expectation is bounded by $1 + c(\lambda)L^{d-1}N^{1-d}$. The estimate (6.31) follows and the proof of Lemma 6.6 is complete. $\square$

We proceed with the proof of Theorem 6.1 when $\alpha \geq \frac{1}{d}$. For any $x \in S_{2N}$ and $m \geq 1$, we find

$$(6.37) \quad \begin{aligned} &P_x^0[U_{B_\infty(0,[N/4])} \leq \mathcal{D}_{[N^\beta]}] \\ &\leq P_x^0[\text{for some } x_* \in B_\infty(0,[N/4]): \tilde{R}_m^{x_*} \leq \mathcal{D}_{[N^\beta]}] \\ &\quad + P_x^0[U_{B_\infty(0,[N/4])} \leq \mathcal{D}_{[N^\beta]}, \\ &\qquad\qquad \text{for all } x_* \in B_\infty(0,[N/4]): \tilde{R}_m^{x_*} > \mathcal{D}_{[N^\beta]}] \\ &\overset{\text{(def.)}}{=} P_1 + P_2. \end{aligned}$$

Applying Lemma 6.6 to $P_1$, we obtain

$$P_1 \quad \leq \quad cN^{d+1} \sup_{x_* \in B_\infty(0,[N/4])} P_x^0[\tilde{R}_m^{x_*} \leq \mathcal{D}_{[N^\beta]}]$$

(6.38)

$$\stackrel{(6.31)}{\leq} \quad cN^{d+1} \exp\{cN^{1-d}L^{d-1}N^\beta - c'm\}.$$

In order to bound $P_2$ in (6.37), we apply the geometric Lemmas 5.1 and 5.3. Together, they imply, for $N \geq c(\gamma,\gamma')$, the following inclusions for the event $A_{\cdot,\cdot,\cdot,\cdot,\cdot}$ defined in (6.13):

$$\{U_{B_\infty(0,[N/4])} \leq \mathcal{D}_{[N^\beta]}, \text{for all } x_* \in S_N : \tilde{R}_m^{x_*} > \mathcal{D}_{[N^\beta]}\}$$

(6.39)

$$\stackrel{\text{(Lemma 5.1)}}{\subseteq} \bigcup_{\substack{x_* \in B_\infty(0,[N/4]), \\ C_{x_*}(L) \subseteq B_\infty(0,[N/4])}} \{X([0,\tilde{D}_m^{x_*}]) \; \tfrac{1}{4}\text{-disconnects } C_{x_*}(L)\}$$

$$\stackrel{\text{(Lemma 5.3)}}{\subseteq} \bigcup_{\substack{x_* \in B_\infty(0,[N/4]), \\ C_{x_*}(L) \subseteq B_\infty(0,[N/4])}} A_{C_{x_*}(L),C_{x_*}(L)^{(L)},l,c(L/l)^d,m}.$$

Since $1 \leq l \leq \frac{2L}{100}$ [cf. (5.1)] and $C_{x_*}(L)^{(2L/10)} \subseteq C_{x_*}(L)^{(L)} \subseteq x_* + S_{2L}$ for $x_* \in B_\infty(0,[N/4])$ and $N \geq c(\gamma,\gamma')$, we can apply Lemma 6.5 with $a = 2L$ to obtain, for $P_2$ in (6.37),

$$P_2 \quad \stackrel{(6.39)}{\leq} \quad N^{d+1} \sup_{x_* \in B_\infty(0,[N/4])} P_x^0[A_{C_{x_*}(L),C_{x_*}(L)^{(L)},l,c(L/l)^d,m}]$$

$$\stackrel{\text{(Lemma 6.5, } a=2L)}{\leq} N^{d+1} \exp\{cm + cL^d l^{-d} \log N - c'L^{d-1}l^{-1}\}.$$

With (6.38) and (6.37), this estimate yields

$$\sup_{x \in S_{2N}} P_x^0[U_{B_\infty(0,[N/4])} \leq \mathcal{D}_{[N^\beta]}]$$

(6.40)

$$\leq cN^{d+1} \exp\{cN^{\beta-(d-1)(1-\gamma')} - c'm\}$$

$$\qquad + N^{d+1} \exp\{cm + cN^{d\gamma'-d\gamma} \log N - c'N^{(d-1)\gamma'-\gamma}\}.$$

In view of (5.1), $0 < \gamma < \gamma' < 1$, and provided $\beta - (d-1)(1-\gamma') < (d-1)\gamma' - \gamma$ and $d\gamma' - d\gamma < (d-1)\gamma' - \gamma$, that is, if

(6.41)          $0 < \gamma < \gamma' < 1, \qquad \beta < d-1-\gamma, \qquad \gamma' < (d-1)\gamma,$

then the right-hand side of (6.40) is bounded from above by $\exp\{-cN^{(d-1)\gamma'-\gamma}\}$ for $m_N = [c''N^{\beta-(d-1)(1-\gamma')}]$ and $N \geq c(\gamma,\gamma')$, for a large enough constant $c'' > 0$. Hence, for $\gamma,\gamma'$ satisfying (6.41), one has, for $N \geq c(\gamma,\gamma')$,

(6.42)          $\displaystyle\sup_{x \in S_{2N}} P_x^0[U_{B_\infty(0,[N/4])} \leq \mathcal{D}_{[N^\beta]}] \leq \exp\{-cN^{(d-1)\gamma'-\gamma}\}.$

For $0 < \beta < d - 1 - \frac{1}{d-1}$, $d \geq 3$, it is easy to check that the constraints (6.41) are satisfied by $\gamma_1 = \frac{1}{d-1} - \frac{\varepsilon_1(d,\beta)}{2(d-1)}$ and $\gamma_1' = 1 - \varepsilon_1(d,\beta)$ for $\varepsilon_1(d,\beta) > 0$ small enough. Moreover, $(d-1)\gamma_1' - \gamma_1 = d - 1 - \frac{1}{d-1} - c\varepsilon_1(d,\beta) = f(\alpha,\beta) - c\varepsilon_1(d,\beta)$. By (6.42), this is enough to show (1.10), since we can make $\varepsilon_1(d,\beta) > 0$ arbitrarily small.

If, on the other hand, $d - 1 - \frac{1}{d-1} \leq \beta < d - 1$, the constraints (6.41) are satisfied by $\gamma_2 = d - 1 - \beta - \frac{\varepsilon_2(d,\beta)}{2(d-1)}$ and $\gamma_2' = (d-1)(d-1-\beta) - \varepsilon_2(d,\beta)$ for $\varepsilon_2(d,\beta) > 0$ sufficiently small and

$$(d-1)\gamma_2' - \gamma_2 = ((d-1)^2 - 1)(d-1-\beta) - c\varepsilon_2(d,\beta)$$
$$\overset{(6.1)}{=} f(\alpha,\beta) - c\varepsilon_2(d,\beta),$$

which yields (1.10) for this range of $\beta$ as well. This completes the proof of Theorem 6.1 for $\alpha \geq \frac{1}{d}$ and hence the proof of Theorem 6.1 altogether. $\quad\square$

REMARK 6.7. It is easy to see from Theorem 1.2 that the exponents in the upper and lower bounds on $T_N^{\mathrm{disc}}$ for $\alpha < 1$ in (1.4) would match if one could show that the large deviations estimate (1.10) holds with the function $f^*$ defined in (1.15). It may therefore be instructive to modify (1.10) by replacing the time $U_{B(\alpha)}$ by $\mathcal{U}$, defined as

$$\mathcal{U} = \inf\{n \geq 0 : X([0,n]) \supseteq \mathbb{T}_N^d \times \{0\}\}.$$

One can then show that $f^*$ is indeed the correct exponent of the corresponding large deviations problem, in the following sense: For any $\alpha, \beta > 0$, $0 < \xi_1 < f^*(\alpha,\beta) < \xi_2$, one has

$$(6.43) \qquad \varlimsup_{N \to \infty} \frac{1}{N^{\xi_1}} \log \sup_{x \in S_{2[N^{d\alpha \wedge 1}]}} P_x^0[\mathcal{U} \leq \mathcal{D}_{[N^{\beta'}]}] < 0 \qquad \text{for any } 0 < \beta' < \beta,$$

as well as

$$(6.44) \qquad \lim_{N \to \infty} \frac{1}{N^{\xi_2}} \log \inf_{x \in S_{2[N^{d\alpha \wedge 1}]}} P_x^0[\mathcal{U} \leq \mathcal{D}_{[N^{\beta'}]}] = 0 \qquad \text{for any } \beta' > \beta.$$

To show (6.43), one notes that standard estimates on one-dimensional random walk imply that the expected amount of time spent by the random walk $X$ in $\mathbb{T}_N^d \times \{0\}$ during one excursion is of order $N^{d\alpha \wedge 1}$. With this information and the observation that $P_\cdot^0[\mathcal{U} \leq \mathcal{D}_{[N^{\beta'}]}] \leq P_\cdot^0[|X([0,\mathcal{D}_{[N^{\beta'}]}]) \cap \mathbb{T}_N^d \times \{0\}| \geq N^d]$, one can apply Khaśminskii's lemma as in the proof of Lemma 6.5 to find the claimed upper bound. For (6.44), one can follow a similar route as in the derivation of the upper bounds on $T_N^{\mathrm{disc}}$. One can first establish Lemma 3.5 and hence the estimate (3.26) with $\infty$ replaced by $\mathcal{D}_1$, and then show that for $\bar{S}_\cdot$ defined in (3.7) and $a_N$ as in (3.32), $P_\cdot^0[\bar{S}_{[c_1 a_N]} \leq \mathcal{D}_1] \geq$ $(1 - cN^{-d\alpha \wedge 1})^{c_1 a_N} \geq c \exp\{-c' N^{d-(d\alpha \wedge 1)}(\log N)^2\}$, where the first inequality follows essentially from standard estimates on one-dimensional random walk. This is enough for (6.44) with $\beta < d - (d\alpha \wedge 1)$. For $\beta \geq d - (d\alpha \wedge 1)$, one uses again that the expected number of visits to $\mathbb{T}_N^d \times \{0\}$ during one excursion is of order $N^{d\alpha \wedge 1}$, and finds that $P^0[\bar{S}_{[c_1 a_N]} < \mathcal{D}_{[N^{\beta'}]}] \to 1$ as $N \to \infty$ for any $\beta' > \beta \geq d - (d\alpha \wedge 1)$. Using (3.26) for the second inequality, one deduces that for $N \geq c(\beta')$,

$$P^0[\mathcal{C}^V_{\mathbb{T}_N^d} > [c_1 a_N] | \bar{S}_{[c_1 a_N]} < \mathcal{D}_{[N^{\beta'}]}] \leq 2P^0[\mathcal{C}^V_{\mathbb{T}_N^d} > [c_1 a_N]] \leq \frac{2}{N^{10}},$$

hence

$$P^0[\mathcal{U} \leq \mathcal{D}_{[N^{\beta'}]}] \geq P^0[\mathcal{C}^V_{\mathbb{T}_N^d} \leq [c_1 a_N] | \bar{S}_{[c_1 a_N]} < \mathcal{D}_{[N^{\beta'}]}] P^0[\bar{S}_{[c_1 a_N]} < \mathcal{D}_{[N^{\beta'}]}]$$
$$\to 1,$$

thus (6.44) for $\beta \geq d - (d\alpha \wedge 1)$. Note that (6.44) and $\{\mathcal{U} \leq \mathcal{D}_{[N^{\beta'}]}\} \subseteq \{U_{B(\alpha)} \leq \mathcal{D}_{[N^{\beta'}]}\}$ together imply that, for any function $f$ in the estimate (1.10), one has $f(\alpha, \beta') \leq f^*(\alpha, \beta)$ for any $\alpha, \beta > 0, \beta' > \beta$, so that $f(\alpha, \beta) \leq f^*(\alpha, \beta)$ whenever $f(\alpha, \cdot)$ is right-continuous at $\beta$.

**Acknowledgments.** The author is indebted to Alain-Sol Sznitman for suggesting the problem and for fruitful advice throughout the completion of this work. Thanks are also due to Laurent Goergen for pertinent remarks on a previous version of this article.

DEPARTMENT MATHEMATIK
ETH ZÜRICH
CH-8092 ZÜRICH
SWITZERLAND
E-MAIL: david.windisch@math.ethz.ch